\documentclass[12pt]{article}[smallcondensed]
\usepackage{pgf,tikz}
\usetikzlibrary{arrows}
\usepackage[USenglish]{babel}
\usepackage{amsmath,amsfonts,amssymb,amsthm}
\usepackage[utf8]{inputenc}
\usepackage[all]{xy}
\usepackage[T1]{fontenc}
\usepackage{xfrac}
\usepackage{authblk}
\usepackage{mathtools}
\usepackage{pifont}
\usepackage{marvosym}

\usepackage{geometry}
\usepackage{tikz}
\usetikzlibrary{matrix}
\usepackage{tikz-cd}
\newtheorem{teo}{Theorem}[section]
\newtheorem{defi}[teo]{Definition}
\newtheorem{prop}[teo]{Proposition}
\newtheorem{lema}[teo]{Lemma}
\newtheorem{coro}[teo]{Corollary}
\newtheorem{obs}{Remark}
\newtheorem*{ejem}{Example}
\theoremstyle{definition}
\newtheorem*{dem}{Proof}

\providecommand{\keywords}[1]
{
  \small	
  \textbf{\textit{Keywords---}} #1
}

\title{\textbf{Information entropy re-defined in a category theory context using preradicals}}

\author[1,3]{\small Sebastian Pardo G}
\author[1,2,3 \Letter]{Gabriel A. Silva}

\affil[1]{Department of Bioengineering\\University of California, San Diego, USA.}
\affil[2]{Department of Neurosciences\\University of California, San Diego, USA.}
\affil[3]{Center for Engineered Natural Intelligence\\University of California, San Diego, USA.}

\date{}

\begin{document}
\maketitle

\renewenvironment{abstract}
{\begin{quote}
\noindent \rule{\linewidth}{.5pt}\par{\bfseries \abstractname.}}
{\medskip\noindent \rule{\linewidth}{.5pt}
\end{quote}
}
\begin{abstract}
Algebraically, entropy can be defined for abelian groups and their endomorphisms, and was latter extended to consider objects in a Flow category derived from abelian categories, such as $R\textit{-}Mod$ with $R$ a ring. Preradicals are endofunctors which can be realized as compatible choice assignments in the category where they are defined. Here we present a formal definition of entropy for preradicals on $R$-Mod and show that the concept of entropy for preradicals respects their order as a big lattice. Also, due to the connection between modules and complete bounded modular lattices, we provide a definition of entropy for lattice preradicals, and show that this notion is equivalent, from a functorial perspective, to the one defined for module preradicals. 
\\
\end{abstract}

Corresponding author: GS, gsilva@ucsd.edu\\



\section{Introduction}

This work is the second of two related papers that aim to use Category Theory as a natural way to describe information flows in networks and relating these ideas to notions of Shannon information theory. In \cite{Flow} we showed that preradicals naturally describe the flow of information. To that end, we used quiver representation theory to model a directed graph as a diagram with objects and morphisms within the category $R$-Mod, for $R$ a principal ideal domain. Particularly, we showed that preradicals generalize the concept of persistence in filtrations, whose underlying structure has the form of a directed acyclic graph.  
 
Broadly, a preradical is an endofunctor \footnote{An endofunctor is a functor with the same domain and codomain} which acts as a subfunctor of the identity functor. Yet, as described in \cite{Sebas}, preradicals can be considered as compatible choice assignments in the category where they are defined. Here, this property motivated us to study the dynamics that  preradicals induce on the objects and their corresponding endomorphisms. To  approach this, we will rely on concepts of entropy defined on abelian groups and their endomorphisms, and on concepts of entropy defined on abelian categories. 

The notion of entropy was first introduced in the field of thermodynamics as a measure of the system's disorder. It was later used in quantum mechanics where it was endowed with a deep mathematical formalism. A third perspective that came later, known as Shannon entropy, measures the efficiency of a system in transmitting information. When it comes to a mathematical subject, entropy was first defined in topology in \cite{Top ent} by Adler, Konheim and McAndrew, for a continuous self-map of a compact space. However, from an algebraic perspective, Weiss settles in \cite{We} the notion of \emph{algebraic entropy} for an endomorphism of an abelian group. On the other hand, in \cite{Pe} Peter provides a slightly different definition of entropy for an automorphism of an abelian group, which was later extended in \cite{Entr Ab} to consider endomorphisms of abelian groups. One substantial difference between Peter's and Weiss's definitions is that, in the former one considers a supremum over all finite subsets of a group $G$ while in the latter one takes the supremum over all finite subgroups. 

Afterwards, algebraic entropy was generalized in \cite{Luigi} to consider not only abelian groups (or $\mathbb{Z}$-modules), but modules over a unitary ring $R$. In this case, one considers \emph{invariants} on $R$-Mod to define an algebraic entropy: invariants are extended real valued functions that are invariant under isomorphic objects. Specially, we focus on invariants $i$ with at least two minimal requirements with which one can associate an algebraic entropy to it. These two conditions comprise what we normally call as \textit{subadditive invariant}. In this work, we will define the notion of entropy for a preradical in terms of an algebraic $i$-entropy, where $i$ represents an invariant defined on $R$-Mod. 

Lastly, due to the relation between lattice preradical and module preradicals described in \cite{Sebas}, we use functors to provide a notion of entropy for lattice preradicals, when restricted to a non-full subcategory $\mathcal{SL_{M}}$ of the category $\mathcal{L_{M}}$ of complete bounded modular lattices and linear morphisms. This functorial connection show that the definition of entropy for lattice preradicals is order-respectful, and further, it is equivalent to the entropy defined for module preradicals. 

The paper is organized as follows: in Section \ref{Section 2} we give some preliminary facts about the structure of the category of flows $\textbf{Flow}_{\mathcal{C}}$ for any category $\mathcal{C}$. We also give a brief introduction to preradicals on $R$-Mod along with notation, terminology, and basic properties. Section \ref{Preradicals in flow} shows that any preradical defined on a category $\mathcal{C}$ induces a preradical on the category of flows associated to $\mathcal{C}$. In particular, we describe the $\alpha$ and $\omega$ preradicals in the category of flows \textbf{Flow}$_{R\textit{-}Mod}$ with respect to the category $R$-Mod. In Section \ref{Section Entropia} we define entropy for preradicals on $R$-Mod and display some examples. Right after, we describe the entropy for preradicals on the category of flows  \textbf{Flow}$_{R\textit{-}Mod}$. Finally, we use some functorial properties between the big lattice of lattice preradicals and the big lattice of module preradicals to define entropy for lattice preradicals, and prove some of its properties.

\section{Preliminaries } \label{Section 2}

\subsection{The Category of Flows}

Given a category $\mathcal{C}$, the identity functor on $\mathcal{C}$ defines the special comma category $(\mathcal{C} \ \downarrow \ \mathcal{C})$. This category, also known as the category of arrows, is usually denoted by $\mathcal{C}^{2}$. The objects of the category $\mathcal{C}^{2}$ are triplets $(C_{1},C_{2},\eta)$, where $C_{1},C_{2}$ are objects of $\mathcal{C}$ and $\eta:C_{1}\longrightarrow C_{2}$ is a morphisms in $\mathcal{C}$. A morphism between two objects $(C_{1},C_{2}, \eta)$ and $(D_{1},D_{2}, \psi)$ of $\mathcal{C}^{2}$ is made of a pair $(f,g)$ where $f:C_{1}\longrightarrow D_{1}$ and $g:C_{2}\longrightarrow D_{2}$ are morphisms in $\mathcal{C}$ such that the diagram 
\begin{equation*}
\xymatrix{
C_{1}\ar[r]^{\eta} \ar[d]_{f} & C_{2} \ar[d]^{g} \\
D_{1} \ar[r]^{\psi}  & D_{2} 
}
\end{equation*}
commutes in $\mathcal{C}$. This arrow category is in fact isomorphic to the functor category $\mathcal{C}^{I}$, where $I$ is the interval category consisting of just two objects and only one nontrivial morphisms connecting them $\{0\longrightarrow 1\}$. 

We can impose two restrictions to the category $\mathcal{C}^{2}$ to obtain a particular subcategory denoted as \textbf{Flow}$_{\mathcal{C}}$. The first restriction consists on taking as triplets  $(C_{1},C_{2},\eta)$ those that satisfies $C_{1}=C_{2}$. This way, the objects of \textbf{Flow}$_{\mathcal{C}}$ are objects in $\mathcal{C}$ equipped by endomorphisms. 
\begin{defi}
A \textit{flow} in $\mathcal{C}$ is a pair $(C, \varphi)$, where $C$ is an object in $\mathcal{C}$ and $\varphi:C\longrightarrow C$ is an endomorphisms in $\mathcal{C}$. 
\end{defi}
For the second restriction, we will take as morphisms those pairs $(f,g)$ satisfying $f=g$. Thus, a \textit{flow morphism} between the flows $(C,\varphi)$ and $(D,\psi)$ is a morphism $f:C\longrightarrow D$ in $\mathcal{C}$ such that the diagram 
\begin{equation*}
\xymatrix{
C\ar[r]^{\varphi} \ar[d]_{f} & C \ar[d]^{f} \\
D \ar[r]^{\psi}  & D
}
\end{equation*}
commutes in $\mathcal{C}$. Therefore,  \textbf{Flow}$_{\mathcal{C}}$ is the category whose objects are flows in $\mathcal{C}$ and whose morphisms are flow morphisms.  

We notice that the above construction  implies that the category \textbf{Flow}$_{\mathcal{C}}$ will not be a full subcategory of $\mathcal{C}^{2}$, although it can inherit some properties of the category $\mathcal{C}$. For instance, as shown in \cite{Entro}, the category \textbf{Flow}$_{\mathcal{C}}$ is isomorphic to the functor category $\mathcal{C}^{\mathbb{N}}$, where $\mathbb{N}$ is considered as a one object category. Thus, the category \textbf{Flow}$_{\mathcal{C}}$ will be abelian when so is $\mathcal{C}$. 

There is a 'forgetful' functor $U:\textbf{Flow}_{\mathcal{C}} \longrightarrow {\mathcal{C}}$ given by the assignment $(C,\eta)\longmapsto C$ on objects and the assignment $f\longmapsto f$ on morphisms. This functor $U$ ignores the dynamic induced by the endomorphisms, since:
\begin{equation*}
\xymatrix{
C\ar[r]^{\eta} \ar[d]_{f}  & C \ar[d]^{f}  & \overset{U}{\longmapsto} & C \ar[d]^{f} \\
D \ar[r]^{\mu}  & D & & D
.} 
\end{equation*}
There is also an embedding functor $E:\mathcal{C}\longrightarrow \textbf{Flow}_{\mathcal{C}}$ given by the assignment $C\longmapsto (C,Id_{C})$ on objects and the assignment $f\longmapsto f$ on morphisms, that is, 
\begin{equation*}
\xymatrix{
C \ar[d]^{f}  & \overset{E}{\longmapsto} & C\ar[r]^{Id} \ar[d]_{f} & C \ar[d]^{f}   \\
D  &  & D\ar[r]^{Id} & D
.} 
\end{equation*}

The next result will be useful for characterizing the subobjects in \textbf{Flow}$_{\mathcal{C}}$; especially when inducing preradicals from $R$-Mod to \textbf{Flow}$_{R\textit{-}Mod}$. 

\begin{lema}\label{mono}
A flow morphism $f$ in \textbf{Flow}$_{\mathcal{C}}$ is a monomorphisms if $U(f)=f$ is a monomorphism in the category $\mathcal{C}$. 
\end{lema}
\begin{dem}
For the sake of the proof, we will denote by $\bar{f}$ the morphism in \textbf{Flow}$_{\mathcal{C}}$ so that $U(\bar{f})=f$. Let $f$ be a monomorphism in $\mathcal{C}$ and let $\bar{f}$ be the flow morphism such that $U(\bar{f})=f$. Suppose that $\bar{f}:(C,\eta)\longrightarrow (D,\mu)$. If $\bar{g}_{1},\bar{g}_{2}:(B,\xi) \longrightarrow (C,\eta)$ are two flow morphisms such that $\bar{f}\circ \bar{g_{1}}=\bar{f}\circ \bar{g_{2}}$, then we have the commutative diagrams 
\begin{equation*} 
\xymatrix{
B\ar[r]^{\xi} \ar[d]_{g_{1}}  & B \ar[d]^{g_{1}}  &  & B \ar[r]^{\xi} \ar[d]_{g_{2}} & B\ar[d]^{g_{2}} \\
C \ar[r]^{\eta} \ar[d]_{f}  & C\ar[d]^{f} & = & C\ar[r]^{\eta} \ar[d]_{f} & C\ar[d]^{f}\\
D\ar[r]^{\mu} & D  &  &  D\ar[r]^{\mu}  & D
} 
\end{equation*}
in $\mathcal{C}$ from which it follows that $(f\circ g_{1} \circ \xi) = (f \circ g_{2} \circ \xi)$. Since $f$ is a monomorphism in $\mathcal{C}$, the last equation simplifies to $g_{1}\circ \xi = g_{2} \circ \xi$. Now, as $\bar{g}_{1},\bar{g}_{2}:(B,\xi) \longrightarrow (C,\eta)$ are both flow morphisms, the condition $g_{1}\circ \xi = g_{2} \circ \xi$ implies that the following diagrams define the same transformation:
\begin{equation*}
\xymatrix{
B\ar[r]^{\xi} \ar[d]_{g_{1}}  & B \ar[d]^{g_{1}}  & \ar @{} [d] |{=} & B\ar[r]^{\xi} \ar[d]_{g_{2}}  & B \ar[d]^{g_{2}} \\
C \ar[r]^{\eta}  & C & & C\ar[r]^{\eta} & C
,} 
\end{equation*}
Therefore, $\bar{g_{1}}=\bar{g_{2}}$.
\qedsymbol
\end{dem}

Formally, a \textit{subobject} of the object $(C,\eta)$ in \textbf{Flow}$_{\mathcal{C}}$, is an equivalence class of monomorphisms into $(C,\eta)$. Thus, we can write a subobject of $(C,\eta)$ as a pair $((X,\alpha),f)$ where $(X,\alpha)\overset{f}{\longrightarrow}(C,\eta)$ is a "representative" monomorphism in \textbf{Flow}$_{\mathcal{C}}$. By Lemma \ref{mono},  it suffices to consider a monomorphism $X\overset{f}{\longrightarrow}C$ in the category $\mathcal{C}$. Hence, for the case $\mathcal{C}=R\textit{-}Mod$, we might choose, as a representative of the equivalence class of $(X,f)$, the pair $(f(X),\iota)$ where $\iota: f(X)\hookrightarrow C$ denotes the inclusion map. Observe that, since $f$ is a monomorphism in the category $R$-Mod, then $f$ is an injective morphism and thus $X\cong f(X)$. Further, by the commutativity of the diagram   
\begin{equation*}
\xymatrix{
X\ar[r]^{\alpha} \ar[d]_{f} & X \ar[d]^{f} \\
C\ar[r]^{\eta}  & C
,}
\end{equation*}
induced by the flow morphism $(X,\alpha)\overset{f}{\longrightarrow}(C,\eta)$, we have that $\eta(f(X))=f(\alpha(X))\subseteq f(X)$; in other words, the endomorphism $\eta$ restricts to $f(X)$. Therefore, the diagram 
\begin{equation*}
\xymatrix{
 f(X) \ar@{^{(}->}[d] \ar[r]^{\eta_{|}}   & f(X) \ar@{_{(}->}[d] \\
 C\ar[r]^{\eta}  & C 
}
\end{equation*}
commutes in $R$-Mod, showing that the pair $((f(X),\eta_{|}),\iota)$ is also a representative subobject of $(C,\eta)$ for the equivalence class containing $((X,\alpha),f)$. 

\begin{obs}\label{U functor}
The forgetful functor $U\colon  \textbf{Flow}_{R\textit{-}Mod}\longrightarrow R\textit{-}Mod$ maps each subobject in \textbf{Flow}$_{R\textit{-}Mod}$ to a subobject in $R$-Mod. In other words, if $((f(X),\eta_{|}),\iota)$ is a representative subobject of $(C,\eta)$ for the class of equivalence containing the flow monomorphism $((X,\alpha),f)$, then $((f(X),\iota)$ is a representative subobject of $C$, for the class of equivalence that contains the monomorphism $X\overset{f}{\longrightarrow} C$ in $R$-Mod.
\end{obs}

\subsection{Preradicals}

Broadly, a \textit{preradical} on $\mathcal{C}$ can be defined in terms of the 2\textit{-}category of all endofunctors\footnote{An endofunctor is a functor $T:\mathcal{C}\longrightarrow \mathcal{C}$ whose domain is equal to its codomain} of $\mathcal{C}$. This 2\textit{-}category coincides with the functor category $\mathcal{C}^\mathcal{C}$, whose objects are all endofunctors $\tau:\mathcal{C}\longrightarrow \mathcal{C}$, and whose morphisms are natural transformations, usually denoted by $\tau_{1} \Longrightarrow \tau_{2}$. Thus, a preradical $\sigma$ on $\mathcal{C}$ is the same as a \textit{co-pointed} endofunctor of $\mathcal{C}^{\mathcal{C}}$; this is, an endofunctor $\sigma:\mathcal{C}\longrightarrow \mathcal{C}$ for which there is natural transformation $\sigma \Longrightarrow Id_{\mathcal{C}}$, and where $Id_{\mathcal{C}}$ denotes the identity functor. Observe that $\sigma \Longrightarrow Id_{\mathcal{C}}$ means that one has a family of morphisms $\{\eta_{C}\}_{C\in \mathcal{C}}$ -indexed by the objects of $\mathcal{C}$- such that, for each $C,D\in \mathcal{C}$ and each morphism  $C\overset{f}{\longrightarrow} D$ in $\mathcal{C}$, the following diagram commutes 
\begin{equation*}
\xymatrix{
C\ar[r]^{f} & D \\
\sigma(C) \ar[r]^{\sigma(f)} \ar[u]^{\eta_{C}} & \sigma(D)\ar[u]_{\eta_{D}}
.}
\end{equation*}

However, according to \cite{Pre 1} one can also interpret a preradical on $\mathcal{C}$ as a subfunctor of the identity functor $Id_{\mathcal{C}}$. In other words, a preradical is a functor that assign to each object $A$ in $\mathcal{C}$, a subobject $\sigma(A)$ in $\mathcal{C}$ such that, for any morphism $A\overset{f}{\longrightarrow} B$ in $\mathcal{C}$, the restriction and corestriction of $f$ to the respective subobjects defines a morphisms $\sigma(f)\colon \sigma(A)\overset{f_{|}}{\longrightarrow} \sigma(B)$ in $\mathcal{C}$, which makes the following diagram to commute  
\begin{equation*}
\xymatrix{
A\ar[r]^{f} & B \\
\sigma(A) \ar[r]^{\sigma(f)} \ar[u]^{\iota_{A}} & \sigma(B)\ar[u]_{\iota_{B}}
.}
\end{equation*}
Notice that in the above diagram, $\iota_{A}$ and $\iota_{B}$ denote the respective inclusion maps. In this way, as the authors describe in  \cite{Sebas}, preradicals behave as compatible choice assignments, since these assign to each object of the category $\mathcal{C}$ a subobject, in such a way that the choices are compatible with respect to the morphism in $\mathcal{C}$. In fact, this property is what makes preradicals useful for describing the flow of information (see \cite{Flow}). 

In this work we will consider preradicals on a category $\mathcal{C}$ as subfunctors of the identity functor $Id_{\mathcal{C}}$. Specially, we will start considering preradicals on $R$-Mod for a unitary ring $R$. Here, objects of $R$-Mod corresponds to $R$-modules and subobjects corresponds to $R$-submodules. For a complete introduction to preradicals and its properties on $R$-Mod see \cite{Bican}, \cite{Pre 1} and \cite{Partitions}. 

\begin{defi}
A preradical $\sigma$ on the category $R$-Mod is a functor $\sigma\colon R\textit{-Mod}\longrightarrow R\textit{-Mod}$ that assigns to each module $M\in R\textit{-Mod}$, a submodule $\sigma(M)$ such that for each morphism $f\colon M\longrightarrow M'$ in $R\textit{-Mod}$, we have the commutative diagram 
\begin{equation*}
\xymatrix{
M\ar[r]^{f} & M' \\
\sigma(M) \ar[r]^{\sigma(f)} \ar[u]^{\iota} & \sigma(M')\ar[u]_{\iota}
.}
\end{equation*}
Here, $\sigma(f)$ is the restriction and corestriction of $f$ to $\sigma(M)$ and $\sigma(M')$ respectively, that is, 
\begin{center}
$\sigma(f)\colon =f\mid_{\sigma(M)}:\sigma(M)\longrightarrow \sigma(M')$.     
\end{center}
Also, $\iota$ represents the inclusion map.
\end{defi}

\begin{ejem}
Consider the category $\mathbb{Z}$-Mod of all $\mathbb{Z}$-modules. This category is isomorphic to the category $Ab$ of all abelian groups. Now, given $M\in \mathbb{Z}$-Mod, one can define   
\begin{center}
$\sigma(M)=\{x\in M \mid 3x=0\}$.
\end{center}

Note that $\sigma(M)$ is a submodule of $M$. Also, for any morphism $f:M\longrightarrow M'$ in $\mathbb{Z}$-Mod and for any $y\in M$, one has that $f(3y)=3f(y)$. Thus, if $x\in \sigma(M)$ then $3f(x)=f(3x)=0$, which implies that $f(x)\in \sigma(M')$.
Hence, $f(\sigma(M))\subseteq \sigma(M')$ which in turns implies that
\begin{equation*}
\xymatrix{
M\ar[r]^{f} & M' \\
\sigma(M) \ar[r]^{f_{|}} \ar[u]^{\iota} & \sigma(M')\ar[u]_{\iota}
}
\end{equation*}
commutes in $\mathbb{Z}\texttt{-}Mod$. 
\end{ejem}

In the category $R$-Mod one can define four principal operations between preradicals: if $M\in R$-Mod and $\sigma,\tau$ are two preradicals, then

\begin{itemize}
\item[i)] $(\sigma \wedge \tau)(M)=\sigma(M)\cap \sigma(M)$;
\item[ii)] $(\sigma \vee \tau)(M)= \sigma(M) + \tau(M)$;
\item[iii)] $(\sigma \cdot \tau)(M)=\sigma (\tau(M))$;
\item[(iv)] $(\sigma:\tau)(M)$ is the submodule of $M$ such that 
\begin{center}
$(\sigma:\tau)(M)/\sigma(M)=\tau(M/\sigma(M))$.
\end{center} 
\end{itemize}

The above operations are called the meet, the join, the product and the coproduct respectively. We note that the product $(\sigma\cdot \tau)$ corresponds to the composition between functors $\tau$ and $\sigma$; while the coproduct $(\sigma:\tau)$ involves taking a quotient module (induced by $\sigma$), apply functor $\tau$, and then use the Correspondence Theorem for Modules to obtain the submodule $(\sigma:\tau)(M)$ of $M$. As regards the join and the meet operations, these can be extended to consider an arbitrary family of preradicals. In other words, for any family of preradicals $\{\tau_{\alpha}\}_{\alpha \in I}$ for $M$ an $R$-module,
\begin{itemize}
\item[i)] $\big( \bigvee_{\alpha\in I} \ \tau_{\alpha} \big)(M) = \Sigma_{\alpha\in I} \ \tau_{\alpha}(M)$,

\item[ii)] $\big( \bigwedge_{\alpha\in I} \ \tau_{\alpha} \big)(M) = \cap_{\alpha\in I} \ \tau_{\alpha}(M)$.
\end{itemize}

We denote by $R$-pr the collection of all preradicals on $R$-Mod. There is a natural partial ordering in $R$-pr given by $\sigma \leq \tau$ if and only if $\sigma(M)\leq \tau(M)$ for each $M\in R$-Mod. This partial ordering together with the meet and join operations make of $R$-pr a big lattice. \footnote{A big lattice is a class (not necessarily a set) having joins and meets for arbitrary families (indexed by a class) of elements}

We close this section with the description of two particular preradicals on $R$-Mod. Let $M$ be an $R$-module and let $N$ be a submodule of $M$. On the one hand, we define the preradical $\alpha_{N}^{M}$ which evaluated in a $K\in R$-Mod is
\begin{center}
$\alpha_{N}^{M}(K)= \sum \{f(N)\mid f\colon M \longrightarrow K \}$.
\end{center}
Observe that $f(N)$ is a submodule of $K$ for each $M\overset{f}{\longrightarrow} K$. On the other hand, we define the preradical $\omega_{N}^{M}$ which evaluated in $K\in R$-Mod is
\begin{center}
$\omega_{N}^{M}(K)= \cap \{f^{-1}(N)\mid f\colon K \longrightarrow M\}$.
\end{center}
Notice here that $f^{-1}(N)$ is a submodule of $K$, for every $f\colon K\longrightarrow M$.

\section{Preradicals in the Category of Flows}\label{Preradicals in flow}

In this section we will show that each preradical on $R$-Mod induces a preradical on the category of flows \textbf{Flow}$_{R\textit{-}Mod}$ and viceversa. We start by noticing that, given a flow $C\overset{\eta}{\longrightarrow} C$ and a preradical $\sigma$ on $R$-Mod, the diagram 
\begin{equation*}
\xymatrix{
C\ar[r]^{\eta} & C \\
\sigma(C) \ar[r]^{\eta|_{\sigma(C)}} \ar@{^{(}->}[u]^{\iota}
& \sigma(C) \ar@{_{(}->}[u]_{\iota}
}
\end{equation*}
commutes in $R$-Mod. This in turn implies that $\eta \mid_{\sigma(C)}:\sigma(C)\longrightarrow \sigma(C)$ is an endomorphism in $R$-Mod, and thus the pair $(\sigma(C),\eta \mid_{\sigma(C)})$ is an object in \textbf{Flow}$_{R\textit{-}Mod}$. Furthermore, we have that $\big(\sigma(C),\eta \mid_{\sigma(C)}\big)$ together with the inclusion map $\iota$ defines a subobject of the flow object $(C,\eta)$ in \textbf{Flow}$_{R\textit{-}Mod}$. 
% Observe that the latter statement is also true for fully invariant\footnote {A submodule $C'$ of $C$ is termed fully invariant if for any endomorphisms $\eta:C\longrightarrow C$ one has $\eta(C')\subseteq C'$.} submodules: if $C'$ is a fully invariant submodule of $C$, then for any endomorphisms $\eta$ of $C$ the pair $(C',\eta |_{C'})$  defines a subobject of the flow object $(C,\eta) \textbf{Flow}_{R\textit{-}Mod}$.

\begin{prop} \label{induced preradicals Flow}
Let $\sigma$ be a preradical on $R$-Mod, and let \textbf{Flow}$_{R\textit{-}Mod}$ be the category of flows associated to $R$-Mod. Then $\sigma$ induces a preradical $\bar{\sigma}$ on \textbf{Flow}$_{R\textit{-}Mod}$ in the following way: $\bar{\sigma}$ assigns to each object $(C,\eta)$ the subobject $(\sigma(C),\eta_{|_{\sigma(C)}})$, and assigns to each flow morphism $(C,\eta)\overset{f}{\longrightarrow} (D,\mu)$ the flow morphism $\sigma(f)$ given by 
\begin{center}
$\sigma(f)\colon (\sigma(C),\eta_{|_{\sigma(C)}}) \longrightarrow  (\sigma(D),\mu_{|_{\sigma(D)}})$   
\end{center}
where $\sigma(f)=f_{|_{\sigma(C)}}$.
\end{prop}

\begin{dem}
Let $\sigma$ be a preradical on $R$-Mod, and let $(C,\eta)$ and $(D,\mu)$ be objects in \textbf{Flow}$_{R\textit{-}Mod}$ with $(C,\eta)\overset{f}{\longrightarrow} (D,\mu)$ a flow morphism. For the sake of this prove, we will denote $(\sigma(C),\eta_{|_{\sigma(C)}})$ by $(\sigma(C),\eta_{|})$ when the context allows no confusion. Hence, we will show that 
\begin{equation*} 
\xymatrix{
(C,\eta)\ar[r]^{f} & (D,\mu) \\
(\sigma(C),\eta_{|}) \ar[r]^{f|_{\sigma(C)}} \ar[u]^{\iota} & (\sigma(D),\mu_{|}) \ar[u]_{\iota}
}
\end{equation*}
is a commutative diagram in \textbf{Flow}$_{\mathcal{C}}$.

As we noticed previously, as $\eta:C\longrightarrow C$ and $\mu:D\longrightarrow D$ are both endomorphisms in $R$-Mod, then we have the following commutative diagrams 
\begin{equation*}
\xymatrix{
C\ar[r]^{\eta} & C \\
\sigma(C) \ar[r]^{\eta|_{\sigma(C)}} \ar[u]^{\iota} & \sigma(C)\ar[u]_{\iota}
}
\mbox{   and    } 
\xymatrix{
D\ar[r]^{\mu} & D \\
\sigma(D) \ar[r]^{\mu|_{\sigma(D)}} \ar[u]^{\iota} & \sigma(D).\ar[u]_{\iota}
}
\end{equation*}
From each of these diagrams we infer that $(\sigma(C),\eta \mid_{\sigma(C)})$ and $(\sigma(D),\mu |_{\sigma(D)})$ are subobjects of $(C,\eta)$ and $(D,\mu)$ respectively. Thus, in order to show that $\sigma$ induces a preradical in \textbf{Flow}$_{\mathcal{C}}$, it suffices to prove that $f|_{\sigma(C)}$ is a flow morphisms from $(\sigma(C),\eta \mid_{\sigma(C)})$ to $(\sigma(D),\mu |_{\sigma(D)})$.

On the one hand, since $\sigma$ is a preradical on $R$-Mod and $f:C\longrightarrow D$ is a morphism in $R$-Mod, then $f|_{\sigma(C)}\colon \sigma(C)\longrightarrow \sigma(D)$ is also a morphism in $R$-Mod. On the other hand, for the flow morphisms $f:(C,\eta)\longrightarrow (D,\mu)$ we have the commutative diagram 
\begin{equation*}
\xymatrix{
C\ar[r]^{\eta} \ar[d]_{f} & C \ar[d]^{f} \\
D \ar[r]^{\mu}  & D 
}
\end{equation*}
from which we obtain that $f \circ \eta = \mu \circ f$. Therefore, when restricted to $\sigma(C)$ we get the following commutative diagram
\begin{equation*}
\xymatrix{
\sigma(C)\ar[r]^{\eta|_{\sigma(C)}} \ar[d]_{f|_{\sigma(C)}} & \sigma(C) \ar[d]^{f|_{\sigma(C)}} \\
\sigma(D) \ar[r]^{\mu|_{\sigma(D)}}  & \sigma(D) 
,}
\end{equation*}
which shows that $f|_{\sigma(C)}$ is a flow morphism between the objects  $(\sigma(C),\eta \mid_{\sigma(C)})$ and $(\sigma(D),\mu |_{\sigma(D)})$ in \textbf{Flow}$_{R\textit{-}Mod}$.
\qed
\end{dem}

We will now describe the induced preradicals on \textbf{Flow}$_{R\textit{-}Mod}$ that correspond to the $\alpha$ and $\omega$ preradicals on $R$-Mod.

\begin{prop}
Let $(M,\eta)$ be an object in \textbf{Flow}$_{R\textit{-}Mod}$ and let $(N,\eta_{|_{N}})$ be a subobject of $(M,\eta)$. Then, for any  $(K,\mu)$ in \textbf{Flow}$_{R\textit{-}Mod}$
\begin{center}
    $\bar{\alpha}_{(N,\eta_{|_{N}})}^{(M,\eta)}(K,\mu) =\big( \underset{f}{\sum}\ f(N), \mu_{|} \big)$
\end{center}
defines a preradical on \textbf{Flow}$_{R\textit{-}Mod}$. \\ Here the sum is taken over all flow morphisms 
$(M,\eta)\overset{f}{\longrightarrow} (K,\mu)$. 
\end{prop}

\begin{dem}
Let $(M,\eta)$ be an object in \textbf{Flow}$_{R\textit{-}Mod}$ and let $(N,\eta_{|_{N}})$ be a subobject of $(M,\eta)$. Before we start with the proof, we will notice that $\big( \underset{f}{\sum}\ f(N), \mu_{|} \big)$ defines a subobject of $(K, \mu)$. For any flow morphism $f\colon (M,\eta)\longrightarrow (K,\mu)$ we have a commutative diagram 
\begin{equation*}
\xymatrix{
M\ar[r]^{\eta} \ar[d]_{f} & M \ar[d]^{f} \\
K \ar[r]^{\mu}  & K 
}
\end{equation*}
in $R$\textit{-}Mod, from which it follows that  
\begin{center}
    $\mu(f(N))=(\mu\circ f)(N)=(f\circ \eta)(N)\subseteq f(N)$.
\end{center}
Considering the above, we have a commutative diagram 
\begin{equation*}
\xymatrix{
f(N)\ar[r]^{\mu_{|}} \ar@{_{(}->}[d]_{\iota} & f(N) \ar@{_{(}->}[d]_{\iota} \\
K \ar[r]^{\mu}  & K 
,}
\end{equation*}
which implies that  $(f(N), \mu_{|_{f(N)}})$ is a subobject of $(K,\mu)$. When considering all flow morphisms $f\colon (M,\eta)\longrightarrow (K,\mu)$, the latter argument can also be apply to $\underset{f}{\sum} \ f(N)$,  as one also has that $\mu\big(\underset{f}{\sum} \ f(N)\big)\subseteq \underset{f}{\sum} \ f(N)$. Therefore, we can say that $\big(\underset{f}{\sum} \ f(N), \mu_{|_{\Sigma}} \big)$ is a subobject of $(K,\mu)$ in \textbf{Flow}$_{R\textit{-}Mod}$.

Now, let $(K,\mu)\overset{h}{\longrightarrow} (K',\mu')$ be a flow morphism. On the one hand, observe that for any flow morphisms $f\colon (M,\eta)\longrightarrow (K,\mu)$, the composition $(h\circ f)\colon (M,\eta)\longrightarrow (K',\mu ')$ is also a flow morphisms. Hence, 
\begin{center}
    $$h\Big(\Big\{\underset{f}{\sum}\ f(N) \mid f\colon M\longrightarrow K \mbox { in } \textbf{Flow}_{R\textit{-}Mod} \Big\}\Big)$$
     $$=\Big\{\underset{f}{\sum}\ (h\circ f)(N) \mid f\colon M\longrightarrow K \mbox { in } \textbf{Flow}_{R\textit{-}Mod} \Big \}$$
$$\subseteq  \Big \{\underset{f'}{\sum}\ f'(N) \mid  f'\colon M\longrightarrow K' \mbox { in } \textbf{Flow}_{R\textit{-}Mod} \Big \}.$$ 
\end{center}

On the other hand, since $h$ is a flow morphisms, we have the commutative diagram
\begin{equation*}
\xymatrix{
K\ar[r]^{\mu} \ar[d]_{h} & K \ar[d]^{h} \\
K' \ar[r]^{\mu'}  & K' 
,}
\end{equation*}
which evaluated in $\underset{f}{\sum}\ f(N)$ induces the commutative diagram
\begin{equation*}
\xymatrix{
\underset{f}{\sum}\ f(N) \ar[r]^{\mu_{|}} \ar[d]_{h_{|}} & \underset{f }{\sum}\ f(N) \ar[d]^{h_{|}} \\
\underset{f'}{\sum}\ f'(N) \ar[r]^{\mu'_{|}}  & \underset{f'}{\sum} \ f'(N) 
}
\end{equation*}
in $R$-Mod. This shows that $h_{|}:\Big( \underset{f}{\sum}\ f(N), \mu_{|} \Big) \longrightarrow \Big( \underset{f }{\sum }\ f'(N), \mu'_{|} \Big)$ is a flow morphism which makes 
\begin{equation*}
\xymatrix{
(K,\mu) \ar[r]^{h}  & (K',\mu')  \\
\bar{\alpha}_{(N,\eta_{|_{N}})}^{(M,\eta)}(K,\mu) \ar@{_{(}->}[u]_{\iota}\ar[r]^{h_{|}}  & \bar{\alpha}_{(N,\eta_{|_{N}})}^{(M,\eta)}(K',\mu') \ar@{_{(}->}[u]^{\iota}
}
\end{equation*}
a commutative diagram in \textbf{Flow}$_{R\textit{-}Mod}$. Therefore, $\bar{\alpha}_{(N,\eta_{|_{N}})}^{(M,\eta)}$ is a preradical on \textbf{Flow}$_{R\textit{-}Mod}$.
\qed
\end{dem}

\begin{prop}
Let $(M,\eta)$ be an object in \textbf{Flow}$_{R\textit{-}Mod}$ and $(N,\eta_{|_{N}})$ a subobject of $(M,\eta)$. Then, for any  $(K,\mu)$ in \textbf{Flow}$_{R\textit{-}Mod}$
\begin{center}
    $\bar{\omega}_{(N,\eta_{|_{N}})}^{(M,\eta)}(K,\mu) =\big( \underset{f }{\cap}\ f^{-1}(N), \mu_{|} \big)$
\end{center}
defines a preradical on \textbf{Flow}$_{R\textit{-}Mod}$. 
\\ Here the intersection is taken over all flow morphisms $f\colon (K,\mu)\longrightarrow (M,\eta)$.
\end{prop}

\begin{dem}
Let $(M,\eta)$ be an object in \textbf{Flow}$_{R\textit{-}Mod}$ and $(N,\eta_{|_{N}})$ a subobject of $(M,\eta)$. We will first notice that, for each flow morphism $f\colon (K,\mu) \longrightarrow (M,\eta)$ the pair $\big(f^{-1}(N), \mu_{|}\big)$ defines a subobject of $(K,\mu)$. For any flow morphism $f\colon (K,\mu)\longrightarrow (M,\eta)$, we get a commutative diagram in $R$-Mod of the form  
\begin{equation*}
\xymatrix{
K\ar[r]^{\mu} \ar[d]_{f} & K \ar[d]^{f} \\
M \ar[r]^{\eta}  & M. 
}
\end{equation*}
If we take into consideration the submodule $f^{-1}(N)\subseteq K$, then we have that 
\begin{align}
\begin{split}\label{omega 1}
    f\big(\mu(f^{-1}(N))\big)=(f\circ \mu)(f^{-1}(N))=(\eta\circ f)(f^{-1}(N)) \\
   =\eta \big(f(f^{-1}(N))\big) \subseteq \eta(N)\subseteq N,
\end{split}
\end{align}
where the last subset relation follows from the fact that $(N,\eta_{|_{N}})$ is an object in \textbf{Flow}$_{R\textit{-}Mod}$. Thus, equation (\ref{omega 1}) implies that $\mu(f^{-1}(N))\subseteq f^{-1}(N)$, and so 
\begin{equation*}
\xymatrix{
f^{-1}(N)\ar[r]^{\mu_{|}} \ar@{_{(}->}[d]_{\iota} & f^{-1}(N) \ar@{_{(}->}[d]^{\iota} \\
K \ar[r]^{\mu}  & K 
}
\end{equation*}
is a commutative diagram in $R$-Mod. This shows that  $(f^{-1}(N), \mu_{|})$ defines a subobject of $(K,\mu)$. Now, when considering all flow morphisms $f\colon (K,\mu)\longrightarrow (M,\eta)$, for the submodule $\underset{f}{\cap} \ f^{-1}(N)$ of $K$ we also have that $\mu\big(\underset{f}{\cap} \ f^{-1}(N) \big)\subseteq \underset{f}{\cap} \ f^{-1}(N)$. Hence, we can say that $\big( \underset{f}{\cap}\ f^{-1}(N), \mu_{|} \big)$ defines a subobject of $(K,\mu)$. 

Let $(K,\mu)\overset{h}{\longrightarrow} (K',\mu')$ be a flow morphism. On the one hand, observe that for any flow morphism $f':(K',\mu')\longrightarrow (M,\eta)$ the composition $(f'\circ h):(K,\mu)\longrightarrow (M,\eta)$ is also a flow morphisms. With this in mind, we have that 
\begin{center}
    $\Big\{\underset{f}{\cap}\ f^{-1}(N) \mid  (K,\mu)\overset{f}{\longrightarrow} (M,\eta) \mbox { in } \textbf{Flow}_{R\textit{-}Mod} \Big\}$ \\ \vspace{2mm}
     $\subseteq \Big\{\underset{f'}{\cap}\ (f'\circ h)^{-1}(N) \mid (K',\mu')\overset{f'}{\longrightarrow} (M,\eta) \mbox { in } \textbf{Flow}_{R\textit{-}Mod} \Big\}$ \\ \vspace{2mm}
$= \Big\{\underset{f'}{\cap}\ h^{-1}(f'^{-1}(N)) \mid  (K',\mu')\overset{f'}{\longrightarrow} (M,\eta) \mbox { in } \textbf{Flow}_{R\textit{-}Mod} \Big\}$ \\ \vspace{2mm}
$= h^{-1}\Big( \Big\{ {\underset{f'}{\cap}\ f'^{-1}(N)  | (K',\mu')\overset{f'}{\longrightarrow} (M,\eta) \mbox { in } \textbf{Flow}_{R\textit{-}Mod} \Big\} \Big)}$,
\end{center}
from which we deduce that 
\begin{center}
    $h\Big(\Big\{\underset{f}{\cap}\ f^{-1}(N) \mid  (K,\mu)\overset{f}{\longrightarrow} (M,\eta) \mbox { in } \textbf{Flow}_{R\textit{-}Mod} \Big\}\Big)$ 
  $\subseteq \Big\{ {\underset{f' \ \ \ \ \ \ \ \ \ \  }{\cap \ f'^{-1}(N)} | (K',\mu')\overset{f'}{\longrightarrow} (M,\eta) \mbox { in } \textbf{Flow}_{R\textit{-}Mod} \Big\} }$.
\end{center}  
On the other hand, since $h$ is a flow morphism, we have the commutative diagram
\begin{equation*}
\xymatrix{
K\ar[r]^{\mu} \ar[d]_{h} & K \ar[d]^{h} \\
K' \ar[r]^{\mu'}  & K' 
}
\end{equation*}
in $R$-Mod. This diagram, when restricted to $\underset{f}{\cap} \ f^{-1}(N)$, induces the commutative diagram
\begin{equation*}
\xymatrix{
\underset{f}{\cap}\ f^{-1}(N) \ar[r]^{\mu_{|}} \ar[d]_{h_{|}} & \underset{f}{\cap}\ f^{-1}(N)  \ar[d]^{h_{|}} \\
\underset{f'}{\cap}\ f'^{-1}(N)  \ar[r]^{\mu'_{|}}  & \underset{f'}{\cap}\ f'^{-1}(N). 
}
\end{equation*}
The latter shows that $h_{|}:\big( \underset{f}{\cap}\ f^{-1}(N) , \mu_{|} \big) \longrightarrow \big( \underset{f'}{\cap}\ f'^{-1}(N) , \mu'_{|} \big)$ is a flow morphism which makes 
\begin{equation*}
\xymatrix{
(K,\mu) \ar[r]^{h}  & (K',\mu')  \\
\bar{\omega}_{(N,\eta_{|_{N}})}^{(M,\eta)}(K,\mu) \ar@{^{(}->}[u]_{\iota}\ar[r]^{h_{|}}  & \bar{\omega}_{(N,\eta_{|_{N}})}^{(M,\eta)}(K',\mu') \ar@{^{(}->}[u]^{\iota}
}
\end{equation*}
a commutative diagram in \textbf{Flow}$_{R\textit{-}Mod}$. Therefore, $\bar{\omega}_{(N,\eta_{|_{N}})}^{(M,\eta)}$ is a preradical on \textbf{Flow}$_{R\textit{-}Mod}$.
\qed
\end{dem}

We will denote by $\textbf{Flow}_{R}\textit{-}pr$ the collection of all preradicals on \textbf{Flow}$_{R\textit{-}Mod}$. As in the case of $R\textit{-}pr$, we can define a partial order on $\textbf{Flow}_{R}$-pr as follows: given $\bar{\sigma}, \bar{\tau}\in \textbf{Flow}_{R}\textit{-}pr$, we say that  $\bar{\sigma}\preceq \bar{\tau}$ if and only if $\bar{\sigma}(M,\eta)$ is a subobject of $\bar{\tau}(M,\eta)$, for each $(M,\eta)\in \textbf{Flow}_{R\textit{-}Mod}$. Observe that by Remark \ref{U functor}, the latter condition amounts to saying that $U\big(\overline{\sigma}(M,\eta)\big)$ is a subobject of $U\big(\overline{\tau}(M,\eta)\big)$ in $R$-Mod. Thus, in this case we might assume that $U\big(\overline{\sigma}(M,\eta)\big)$ is a submodule of $U\big(\overline{\tau}(M,\eta)\big)$. 
\\ As we next see, the assignment defined on Proposition \ref{induced preradicals Flow} is injective and order preserving.

\begin{prop}\label{R-pr and Flow}
Let $\phi:R\textit{-}pr \longrightarrow \textbf{Flow}_{R}\textit{-}pr$ be the assignment defined by  $\sigma \longmapsto \overline{\sigma}$ as in Proposition \ref{induced preradicals Flow}. Then $\phi$ is order preserving and an object-injective assignment.
\end{prop}

\begin{dem}
Let $\sigma,\tau \in R\textit{-}pr$ with $\sigma\leq \tau$, and let $\phi(\sigma)=\overline{\sigma}$ and $\phi(\tau)=\overline{\tau}$ be their respective preradicals on $\textbf{Flow}_{R}\textit{-}pr$. As $\sigma\leq \tau$ in $R\textit{-}pr$, then $\sigma(M)\leq \tau(M)$ for every $R$-module $M$. Moreover, if for any endomorphism $\eta\colon M\longrightarrow M$ we denote by $\overline{\sigma}(M,\eta)=(\sigma(M),\eta_{|})$ and $\overline{\tau}(M,\eta)=(\tau(M),\eta_{|})$, then the  commutative diagram  
\begin{equation}
\xymatrix{
\sigma(M) \ar[r]^{\eta_{|}} \ar@{^{(}->}[d]_{\iota} & \sigma(M) \ar@{^{(}->}[d]^{\iota} \\
\tau(M) \ar[r]^{\eta_{|}}  & \tau(M) 
}
\end{equation}
in $R\textit{-}Mod$ implies that $(\tau(M),\eta_{|})$ is a subobject of $(\sigma(M),\eta_{|})$ in $\textbf{Flow}_{R}\textit{-}pr$. Thus, $\overline{\sigma}\preceq \overline{\tau}$. 

For the second statement, suppose that $\sigma,\tau\in R\textit{-}pr$ are such that 
\begin{center}
$\overline{\sigma}=\phi(\sigma)=\phi(\tau)=\overline{\tau}$.    
\end{center}
Then for any object $(M,\eta)\in \textbf{Flow}_{R\textit{-}Mod}$ we have that 
\begin{center}
$(\sigma(M),\eta_{|})=\overline{\sigma}(M,\eta)=\overline{\tau}(M,\eta)=(\tau(M),\eta_{|})$,    
\end{center}
which in turn implies that $\sigma(M)=\tau(M)$. Since this happens for any $(M,\eta)\in \textbf{Flow}_{R\textit{-}Mod}$, we can infer that $\sigma(M)=\tau(M)$ for all $M\in R\textit{-}Mod$. Therefore, $\sigma=\tau$ in $R\textit{-}pr$, and thus, $\phi$ is an object-injective assignment.
\qed
\end{dem}

Next we show that every preradical on $\textbf{Flow}_{R}\textit{-}pr$ induces a preradical in $R\textit{-}pr$. This induced preradical is obtained via the forgetful functor $U$ together with the embedding functor $E$. Before we begin with the proof, note that for any $M\in R$-Mod one can consider the object $E(M)=(M,Id_{M})$ in \textbf{Flow}$_{R\textit{-}Mod}$, and then evaluate any preradical $\overline{\sigma}\in \textbf{Flow}_{R}\textit{-}pr$ on $(M,Id_{M})$. This way we will get a subobject $\overline{\sigma}(M,Id_{M})=(M',Id_{M'})$ of $(M,Id_{M})$. Also observe that, by the way subobjects are defined in  \textbf{Flow}$_{R\textit{-}Mod}$, we can assume that $M'\overset{\iota}{\hookrightarrow} M$, with which the following diagram  
\begin{equation}
\xymatrix{
M' \ar[r]^{Id_{M'}} \ar@{^{(}->}[d]_{\iota} & M' \ar@{^{(}->}[d]^{\iota} \\
M \ar[r]^{Id_{M}}  & M 
}
\end{equation}
commutes in $R$-Mod. 

\begin{prop}
Each preradical $\bar{\sigma}$ on \textbf{Flow}$_{R\textit{-}Mod}$ induces a preradical $\sigma$ on $R$-Mod. This preradical is given by $\sigma:=(U\circ \overline{\sigma} \circ E)$ and satisfies that, for each $M\in R$-Mod, 
\begin{center}
    $\sigma(M)=M'$ where $M'$ is such that $\bar{\sigma}(M,Id_{M})=(M',Id_{M'})$.
\end{center}
\end{prop}

\begin{dem}
Let $\bar{\sigma}$ be a preradical on \textbf{Flow}$_{R\textit{-}Mod}$ and let us take $M, L\in R$-Mod with $g\colon M\longrightarrow L$ an $R$-morphism. For the objects $M$ and $L$ in $R$-Mod along with the morphism $g\colon M\longrightarrow L$, their image under the embedding functor $E:R\textit{-}Mod \longrightarrow \textbf{Flow}_{R\textit{-}Mod}$ is given by the objects $(M,Id_{M})$ and $(L,Id_{L})$ together with the flow morphism $g\colon (M,Id_{M})\longrightarrow (L,Id_{L})$. Thus, when applying the preradical $\bar{\sigma}\in  \textbf{Flow}_{R}\textit{-}pr$, we obtain the commutative diagram 
\begin{equation*}\label{induced preradical}
\xymatrix{
(M,Id_{M})\ar[r]^{g}  & (L,Id_{L})  \\
(M',Id_{M'}) \ar[r]^{g_{|}} \ar@{^{(}->}[u]^{\iota} & (L',Id_{L'}) \ar@{^{(}->}[u]_{\iota} 
}
\tag{*}
\end{equation*}
in \textbf{Flow}$_{R\textit{-}Mod}$. Here $\bar{\sigma}(M,Id_{M})=(M',Id_{M'})$ and $\bar{\sigma}(L,Id_{L})=(L',Id_{L'})$. 
Also, observe that  $(M',Id_{M'})\overset{g_{|}}{\longrightarrow} (L',Id_{L'})$ is a flow morphism, which implies that $g_{|}:M'\longrightarrow L'$ is a morphisms in $R$-Mod that makes 
\begin{equation*}
\xymatrix{
M'\ar[r]^{Id_{M'}} \ar[d]_{g_{|}} & M' \ar[d]^{g_{|}} \\
L' \ar[r]^{Id_{L'}}  & L' 
}
\end{equation*}
a commutative diagram. Finally, as all functors preserve commutative diagram, if we apply functor $U$ to diagram (\ref{induced preradical}) we obtain the commutative diagram in $R$-Mod:
\begin{equation*}
\xymatrix{
M\ar[r]^{g}  & L  \\
M' \ar[r]^{g_{|}} \ar@{^{(}->}[u]^{\iota} & L'\ar@{^{(}->}[u]_{\iota} .
}
\end{equation*}
Therefore, $\sigma\colon =(U\circ \overline{\sigma} \circ E)$ defined on objects as $\sigma(M)=M'$ and on morphisms as $\sigma(g)=g_{|}$, defines a preradical on $R$-Mod.
\qed
\end{dem}

\begin{coro}
The assignment $\textbf{Flow}_{R}\textit{-}pr\longrightarrow \textbf{R}\textit{-}pr$ given by $\overline{\sigma}\longmapsto \sigma$ is order preserving. 
\end{coro}
\begin{dem}
Let $\overline{\sigma}$ and $\overline{\tau}$ be preradicals on \textbf{Flow}$_{R\textit{-}Mod}$ such that $\overline{\sigma}\preceq \overline{\tau}$. Thus, for every $M\in R\textit{-}Mod$ one has that $\overline{\sigma}(M,Id_{M})$ is a subobject of $\overline{\tau}(M,
Id_{M})$. As noted before, this is equivalent to saying that $U\big (\overline{\sigma}(M,Id_{M})\big)$ is a subobject of $U\big(\overline{\tau}(M,
Id_{M})\big)$ in $R$-Mod, which in turn is equivalent to $U\big(\overline{\sigma}(M,Id_{M})\big) \leq U\big(\overline{\tau}(M,
Id_{M})\big)$. Therefore, as $E(M)=(M,Id_{M})$ we have that 
\begin{center}
    $\sigma(M)=(U\circ \overline{\sigma} \circ E)(M)=U\big(\overline{\sigma}(M,Id_{M})\big) \leq U\big(\overline{\tau}(M,
Id_{M})\big)=(U\circ \overline{\tau} \circ E)(M)=\tau(M)$.
\end{center}
Thus, $\sigma(M)\leq \tau(M)$ for all $M\in R\textit{-}Mod$, and therefore $\sigma\leq \tau$ in $R\textit{-}pr$. 
\qed
\end{dem}

%We recall from \cite{Entro} that any algebraic flow over an $R$-module $M$, defines a particular $R[t]$-module structure on $M$. On the other hand, any $R[t]$-module $M$ can be realized as the algebraic flow $(M,\mu)$, where $\mu:M\longrightarrow M$ is given by $\mu(m)=tm$. In fact, one has 
%\begin{prop}\cite[Theorem 3.3]{Entro}  
%Let $R$ be a ring. Then $\textbf{Flow}_{R\textit{-}Mod}$ and $R[t]$-Mod are isomorphic categories.
%\end{prop}

%The above proposition implies that each preradical on $\textbf{Flow}_{R\textit{-}Mod}$ corresponds to a preradical on $R[t]$-Mod and viceversa. From Proposition \ref{R-pr and Flow} it follows the next 
%\begin{coro}
%There is an injective order preserving assignment $\phi:R\textit{-}pr \longrightarrow R[[t]\textit{-}pr$.
%\end{coro}

\section{Entropy for preradicals}\label{Section Entropia}

We start this section recalling the notion of an invariant function on $R$-Mod along with some notation and basic facts which lead to the definition of entropy for endomorphisms and objects in $R$-Mod. Afterwards, we will define entropy for preradicals on $R$-Mod after for preradicals on \textbf{Flow}$_{R}\textit{-}pr$ as a consequence. For a complete introduction to invariants and entropy functions on $R$-Mod, see \cite{Intro Alg}, \cite{Entr Ab},\cite{Alg Entr},\cite{Entro}, \cite{Norm} and \cite{Luigi}.

\begin{defi} \cite{Entro}[Definition 1.1] An \textit{invariant} in $R$-Mod is a function $i:R\textit{-}Mod\longrightarrow \mathbb{R}_{\geq 0}\cup \infty$ such that $i(0)=0$ and $i(M)=i(M')$ whenever $M$ and $M'$ are isomorphic objects in $R$-Mod.
\end{defi}

From now on we will consider invariants with the two following conditions:
\begin{itemize}
    \item[(a)] $i(M_{1}+M_{2})\leq i(M_{1})+i(M_{2})$ for all submodules $M_{1},M_{2}$ of $M$.
    \item[(b)] $i(M/N)\leq i(M)$ for any submodule $N$ of $M$.
\end{itemize}
An invariant $i$ satisfying $(a)$ and $(b)$ is called a \textit{subadditive invariant}.

Let $(M,\eta)$ be an object in $\textbf{Flow}_{R\textit{-}Mod}$ and let $L$ be a subset of $M$. For any positive integer $n$, we define the $n$-\textit{th} $\eta$-\textit{trajectory} of $L$ as 
\begin{center}
    $T_{n}(L,\eta)= L+\eta(L)+\cdots + \eta^{n-1}(L)$.
\end{center}
When $L$ happens to be a submodule of $M$, then $T_{n}(L,\eta)$ is in fact a submodule of $M$. Furthermore, the sum \begin{center}
$\underset{n\geq 1}{\sum}  \ T_{n}(L,\eta)=\underset{n\geq 0}{\sum} \ \eta^{n}(L)$    
\end{center}
defines the submodule of $M$ named the $\eta$-\textit{trajectory} of $L$, which is denoted by $T(L,\eta)$. 

Let us now suppose that $i$ is a subadditive invariant on $R$-Mod. Then, for every submodule $L$ of $M$, with $i(L)<\infty$, one has that $i\Big( T_{n}(L,\eta) \Big)<\infty$ for every $n\geq 1$. Indeed, by proceeding with induction over $n$, we will start considering the base case $n=2$ since $n=1$ is trivial. As $\eta(L)\cong \Big(L/ker(\eta)\Big)$ and $i$ is a subadditive invariant, we have that 
\begin{center}
$i\big( \eta(L) \big)=i\Big(L/ker(\eta)\Big)\leq i(L)<\infty$.    
\end{center}
Hence, 
\begin{center}
$i\big( T_{2}(L,\eta) \big)=i\big( L+\eta(L) \big)\leq i(L)+\eta(L)<\infty$.
\end{center}
Assume now that the inductive hypothesis is valid for $n>2$. Thus, since $\eta^{n}(L)\cong \Big( L/ker(\eta^{n}) \Big)$ we have that 
\begin{center}
$i\big( \eta^{n}(L)\big)=i\Big( L/ker(\eta^{n}) \Big)\leq i(L)< \infty$.
\end{center} 
Therefore, as   
\begin{center}
$T_{n+1}(L,\eta)=L+\eta(L)+\cdots + \eta^{n-1}(L)+\eta^{n}(L)=T_{n}(L,\eta)+\eta^{n}(L)$,    
\end{center} 
by the inductive hypothesis we have that
\begin{center}
$i\big(T_{n+1}(L,\eta)\big) = i\Big( T_{n}(L,\eta)+\eta^{n}(L)\Big)$
$\leq i\Big(T_{n}(L,\eta)\Big) + i\big(\eta^{n}(L)\big) < \infty$.
\end{center}

We now observe that, for every subadditive invariant $i$ and any submodule $L$ of $C$, one has
\begin{gather*}
T_{n+m}(C,\eta)=L+\eta(L)+\cdots + \eta^{n-1}(L)+\eta^{n}(L)+\cdots + \eta^{n+m-1}(L) \\
=T_{n}(L,\eta)+ \eta^{n}(L)+\cdots + \eta^{n+m-1}(L) \\ =T_{n}(L,\eta)+\eta^{n}\big( L+\eta(L)+\cdots + \eta^{m-1}(L)\big) 
=T_{n}(L,\eta)+\eta^{n}\big( T_{m}(L,\eta) \big).            
\end{gather*}
Since $i$ is a subadditive invariant, it follows that 
\begin{gather*}
i\Big( T_{n+m}(L,\eta) \Big) = i\Big( T_{n}(L,\eta) + \eta^{n} \big(T_{m}(L,\eta)\big) \Big) \\
\leq i\Big( T_{n}(L,\eta) \Big) + i\Big( \eta^{n}\big(T_{m}(L,\eta)\big)\Big)\\ 
\leq i\Big( T_{n}(L,\eta) \Big) + i\Big( T_{m}(L,\eta) \Big),
\end{gather*}
where the last relation follows from $(b)$ and the fact that $\eta^{n}\big(T_{m}(L,\eta)\big)\cong \Big( T_{m}(L,\eta) / ker(\eta^{n}) \Big)$.
\begin{obs}
If $L$ is a submodule of $M$ such that $i(L)<\infty$, then every endomorphism $\eta$ of $M$ induces a sequence of positive reals $\{a_{i}\}_{n\in \mathbb{N}}$ such that $a_{n+m}\leq a_{n}+a_{m}$.
\end{obs}

We lastly note that every sequence of positive reals $\{a_{i}\}_{n\in \mathbb{N}}$ such that $a_{n+m}\leq a_{n}+a_{m}$ converges in $\mathbb{R}$. Indeed, as $a_{2}\leq a_{1}+a_{1}=2\cdot a_{1}$ by mathematical induction we have that $a_{k}\leq k\cdot a_{1}$ for all $k\geq 1$. This in turn implies that the sequence $\{\frac{a_{k}}{k} | k\geq 1 \}$ is bounded above by $a_{1}$ and bounded below by $0$. Therefore, $inf \{\frac{a_{k}}{k} | k\geq 1 \}$ exist. Furthermore, one exactly has $\underset{n\to \infty}{lim} \ \frac{a_{n}}{n}$ coincides with $inf \{\frac{a_{k}}{k} | k\geq 1 \}$ as the next Proposition states:

\begin{prop}\cite[Exercise 6.5]{Wa}
Let $\{a_{i}\}_{i\in \mathbb{N}}$ be a sequence of positive real numbers such that $a_{n+m}\leq a_{n}+a_{m}$ for all $n,m\in \mathbb{N}$. Then the sequence $\{\frac{a_{k}}{k} | k\geq 1 \}$ converges to $inf \{\frac{a_{k}}{k} | k\geq 1 \}$.
\end{prop}

The above proposition is the prelude that allows giving a formal definition of $i$-entropy for each invariant $i$ on $R$-Mod. In doing this, one must consider the collection of submodules $M'$ of $M$ such that $i(M')<\infty$, which is denoted by 
\begin{center}
    $\mathcal{F}_{i}(M)\colon =\{M' \leq M \mid i(M')<\infty \}$.
\end{center}

\begin{defi}\cite[Definition 5.25]{Norm}
Let $i$ be a subadditive invariant, $M$ an $R$-module and $M\overset{\eta}{\longrightarrow} M$ an endomorphism. The algebraic $i$-entropy of $\eta$ with respect to $M'\in \mathcal{F}_{i}(M)$ is 
\[
H_{i}(M',\eta)\colon = \underset{n\to \infty}{Lim} \frac{i(T_{n}(M',\eta))}{n}.\]
The \textbf{algebraic $i$-entropy} of $\eta$ is 
\begin{center}
$ent_{i}(\eta)\colon =sup\Big \{H_{i}(M',\eta) \mid M'\in \mathcal{F}_{i}(M)\Big \}$.
\end{center}
\end{defi}

\begin{defi}\cite[Definition 3]{Luigi}
Let $i$ be a subadditive invariant of $R\textit{-}Mod$ and let $End(M)$ denote the set of all endomorphism of the module $M\in R\textit{-}Mod$. The \textbf{$i$-entropy of} $M$ is 
\begin{center}
    $ent_{i}(M)\colon =sup \Big \{ent_{i}(\eta) \mid \eta\in End(M) \Big \}
    = \underset{\eta\in End(M)}{sup}\Big \{ent_{i}(\eta) \Big \}$.
\end{center}
\end{defi}

We now present the definition of $i$-entropy for preradicals on $R$-Mod:

\begin{defi}
Let $i$ be a subadditive invariant on $R\textit{-}Mod$, and let $\sigma\in R\textit{-}pr$. If $M\in R\textit{-}Mod$, then the algebraic $i$-entropy of $\sigma$ with respect to $M$ is given by 
\[
ent_{i}(\sigma)|_{M}\colon =sup \Big \{ H_{i}((M',\eta)  \mid M'\in \mathcal{F}_{i}(\sigma(M)) \mbox{ and }\eta\in End(M) \Big \}
\]
\[
=\underset{\eta \in End(M)}{sup} \Big \{H_{i}(M',\eta) \mid M'\in \mathcal{F}_{i}(\sigma(M))\Big \}
\]
where $H_{i}(M',\eta)=\lim_{n\to \infty} \frac{i(T_{n}(M',\eta))}{n}$.
\end{defi}

\begin{obs}
If $i$ is an invariant on $R$-Mod, then for every preradical $\sigma\in R\textit{-pr}$ and for every $M\in R$-Mod one has that 
\[ 
ent_{i}(\sigma)|_{M} \leq_{\mathbb{R}} ent_{i}(\sigma(M)).
\]
\end{obs}

\begin{prop}
Let $i$ be an invariant on $R$-Mod. If $\sigma, \tau\in R\textit{-}pr$ are such that $\sigma\leq \tau$, then for any $M\in R\textit{-}Mod$, one has 
\[
ent_{i}(\sigma)|_{M}\leq_{\mathbb{R}} ent_{i}(\tau)|_{M}.
\]
\end{prop}

\begin{dem}
Let $M\in R\textit{-}Mod$. As $\sigma(M)\leq \tau(M)$, then for any endomorphism $M\overset{\eta}{\longrightarrow}M$ one has that $\mathcal{F}_{i}(\sigma(M))\subseteq \mathcal{F}_{i}(\tau(M))$. Hence,  
\[
\big \{ H_{i}(M',\eta) \ | \  M'\in \mathcal{F}_{i}(\sigma(M)) \big \} \leq_{\mathbb{R}} \big \{ H_{i}(M',\eta) \ | \ M'\in \mathcal{F}_{i}(\tau(M)) \big \}, 
\]
and thus 
\[
\underset{\eta \in End(M)}{sup} \big \{ H_{i}(M',\eta) \ | \  M'\in \mathcal{F}_{i}(\sigma(M)) \big \} \leq_{\mathbb{R}} \underset{\eta \in End(M)}{sup} \big \{ H_{i}(M',\eta) \ | \ M'\in \mathcal{F}_{i}(\tau(M)) \big \}. 
\]
Therefore, $ent_{i}(\sigma)|_{M}\leq_{\mathbb{R}} ent_{i}(\tau)|_{M}$.
\end{dem}

A well known property in the big lattice $R$-pr that relates the four basic operations between preradicals is: given $\sigma,\tau\in R$-pr one always has that 
\begin{center}
$(\sigma \cdot \tau)\leq (\sigma \wedge \tau) \leq (\sigma \vee \tau) \leq (\sigma: \tau)$. 
\end{center}

As a straightforward result, we have the next

\begin{prop}\label{Inequality}
Let $\sigma,\tau \in R\textit{-}pr$ and $M\in R\textit{-}Mod$. Then  
\begin{center}
    $ent_{i}(\sigma\cdot \tau)|_{M}\leq_{\mathbb{R}} 
    ent_{i}(\sigma\wedge \tau)|_{M} \leq_{\mathbb{R}}
    ent_{i}(\sigma\vee \tau)|_{M} \leq_{\mathbb{R}} 
    ent_{i}(\sigma : \tau)|_{M}$.
\end{center}
\end{prop}

Next we display two examples of invariants for a fix preradical on $R$-Mod. 

\begin{ejem}
Let $R=\mathbb{Z}$ -so the category $R$-Mod coincides with the category of abelian groups $\mathbf{Ab}$- and let us consider the preradical $Tor$ on $\mathbb{Z}\text{-Mod}$ which assigns to each abelian group $M$ its torsion subgroup $Tor(M)$, this is, its subgroup consisting of all elements of finite order. 

On the other hand, we can consider the algebraic entropy on $\mathbb{Z}$-Mod defined by the invariant $i\colon \mathbb{Z}\textit{-}Mod\longrightarrow \mathbb{Z}\textit{-}Mod$, where $i(M)=|log(M)|$ whenever $M$ is finite, otherwise $i(M)=\infty$. We recall that if $\mathcal{F}(M)$ denotes all finite subgroups of the abelian group $M$, then for any endomorphism $M\overset{\eta}{\longrightarrow} M$ and any $L\in \mathcal{F}(M)$, we define for each $n\in \mathbb{N}$
\begin{center}
$H_{n}(L,\eta)=log |T_{n}(L,\eta)|$.
\end{center}
In other words, $H_{n}(L,\eta)$ is the logarithm of the $n\textit{-}th$ trajectory. Thus, the algebraic entropy of $\eta$ with respect to the finite subgroup $L$ is
\begin{center}
$H(L,\eta)= \lim_{n\to \infty} \frac{log|(T_{n}(L,\eta))|}{n}$,   
\end{center}
and hence, the algebraic entropy of $\eta$ is
\begin{center}
    $ent_{i}(\eta)=sup\Big \{H(L,\eta) \mid L\in \mathcal{F}(M)\Big \}$.
\end{center}
Taking the latter into consideration, we define the algebraic entropy of $M$ as 
\begin{center}
$ent_{i}(M)=sup \Big \{ ent_{i}(\eta) | \eta\in End(M) \Big \}$.    
\end{center}

Let us define $M=\oplus_{i\in \mathbb{N}} \ \mathbb{Z}_{p}^{i}$, where $\mathbb{Z}_{p}^{i}$ denotes the cyclic group of order $p$, for each $i\in \mathbb{N}$. Note that, as $M$ is clearly a torsion group, we have that $Tor(M)=M$. Furthermore, if $\beta:M\longrightarrow M$, defined by
\begin{center}
$\beta(x_{1},x_{2},x_{3},\cdots )\longmapsto (0,x_{1},x_{2},x_{3},\cdots )$
\end{center}
is the \textbf{Bernoulli shift endomorphism}, then one has that $ent(\beta)=|log(p)|$ and $ent(Tor)\mid _{M} > 0$ as a consequence. Indeed, if we take the subgroup $\mathbb{Z}_{p}^{1}$ of $M$, for each $n\in \mathbb{N}$ we have that
\begin{center}
$H_{n}(\mathbb{Z}_{p}^{1},\beta)=log|\oplus_{i\leq n} \ \mathbb{Z}_{p}^{i}|=log(|p^{n}|)=n\times log(p)$.    
\end{center}
Hence,  
\begin{center}
$\frac{H(\mathbb{Z}_{p}^{1},\beta)_{n}}{n}=log(p)$    
\end{center}
for all $n\in \mathbb{N}$. Thus, $H(\mathbb{Z}_{p}^{1},\beta)=log(p)$, which in turn implies that $ent_{i}(\beta)\geq log(p)$. Moreover, since the $\beta$-trajectory of the finite subgroup $\mathbb{Z}_{p}^{1}$ covers $M$, then it also covers any other $\beta$-trajectory for any finite subgroup $F$ of $M$. Consequently, we have that $H(F,\beta)\leq H(\mathbb{Z}_{p}^{1},\beta)$, and therefore  
\begin{center}
$ent_{i}(\beta)=H_{i}(\mathbb{Z}_{p}^{1},\beta)=log(p)$.    
\end{center}
Observe that the latter implies that $ent_{i}(Tor)\mid _{M}> 0$. 

On the other hand, if we consider the \textbf{rank invariant} on $\mathbb{Z}$-Mod we will get that $ent_{\textbf{rk}}(Tor)|_{M}=0$ for all $M\in \mathbb{Z}\textit{-}$Mod. Recall that the rank invariant is defined as follows: for each $M\in \mathbb{Z}\textit{-}Mod$, $\textbf{rk}(M)=dim_{\mathbb{Q}}(M\otimes \mathbb{Q})$ when this dimension is finite, otherwise we set $\textbf{rk}(M)=\infty$. Thus, we have that $\textbf{rk}(M)=0$ for any torsion group $M$. This in turn implies that 
\begin{center}
$ent_{\textbf{rk}}(Tor)|_{M}=0$     
\end{center}
for all $M\in \mathbb{Z}\textit{-}$Mod.
\end{ejem}

% finally, as it is noticed in \cite{Alg Entr}, the torsion group $\mathbb{Z}_{p}^{\infty}$, with generators $\langle p_{n} \rangle$ for each $n\in \mathbb{N}$, satisfies that $ent_{a}(\mathbb{Z}_{p}^{\infty})=\infty$. The latter follows from the fact that the endomorphism $\alpha:\mathbb{Z}_{p}^{\infty} \longrightarrow \mathbb{Z}_{p}^{\infty}$ defined by $\alpha(\langle p_{k}\rangle)=p\cdot \langle p_{k+1} \rangle$ has $ent_{a}(\alpha)=\infty$. 

% The invariants from above satisfy an extra condition, making them additive invariants on $R$-Mod. These additive properties are compatible with preradicals, for example in the additive Theorem for abelian groups. 

We will now define the $i$-entropy for preradicals on \textbf{Flow}$_{R\textit{-}Mod}$ object-wise, this is, in terms of flows, and then we will prove some correlations between the entropy of a preradical on $R$-Mod with the entropy of its corresponding preradical on the category of flows \textbf{Flow}$_{R\textit{-}Mod}$.

\begin{defi}
Let $i$ be a subadditive invariant on $R\textit{-}Mod$ and let $\overline{\sigma}$ be a preradical on $\textbf{Flow}_{R}\textit{-}pr$. By denoting   $\overline{\sigma}(M,\eta)=(\sigma(M),\eta_{|})$ for each flow object $(M,\eta)\in \textbf{Flow}_{R\textit{-}Mod}$, we have that the algebraic $i$-entropy of $\overline{\sigma}$ with respect to $\eta$ and $M'\in \mathcal{F}_{i}(\sigma(M))$ is 
\[
H_{i}(M',\eta)\colon = \lim_{n\to \infty} \frac{i(T_{n}(M',\eta))}{n}.
\]
The \textbf{algebraic $i$-entropy} of $\overline{\sigma}$ regarding to the object $(M,\eta)\in \textbf{Flow}_{R\textit{-}Mod}$ is 
\begin{center}
    $ent_{i}(\overline{\sigma})_{(M,\eta)}\colon =sup\Big \{H_{i}(M',\eta) \mid M'\in \mathcal{F}_{i}(\sigma(M))\Big \}$.
\end{center}
\end{defi}
In other words, the algebraic $i$-entropy of $\overline{\sigma}$, with respect to the flow object $(M,\eta)$, is the same as the algebraic $i$-entropy of the endomorphism obtained by restricting $\eta$ to $\sigma(M)$, this is,  $\eta_{|}:\sigma(M)\longrightarrow \sigma(M)$. By how $\psi:R\textit{-}pr \longrightarrow \textbf{Flow}_{R}\textit{-}pr$ is defined, we have the straightforward result:

\begin{prop}
Let $\sigma, \tau \in R\textit{-}pr$ and $(M,\eta) \in \textbf{Flow}_{R\textit{-}Mod}$. If $\overline{\sigma}$ and $\bar{\tau}$ denote the preradicals induced by $\sigma$ and $\tau$, respectively, via the assignment $\psi:R\textit{-}pr \longrightarrow \textbf{Flow}_{R}\textit{-}pr$, then for each invariant $i$ on $R$-Mod
\begin{itemize}
    \item The $i$-entropy for $(\overline{\sigma \wedge \tau})$ with respect to $(M,\eta)$ is 
$$ent_{i}(\overline{\sigma \wedge \tau})_{(M,\eta)}=sup\Big\{H_{i}(M',\eta) \mid M'\in \mathcal{F}_{i}\big( (\sigma \wedge \tau) (M)\big) \Big \},$$
that is, where $M'\leq \sigma(M)\cap \tau(M)$ and $i(M')<\infty$.

\item The $i$-entropy for $(\overline{\sigma \vee \tau})$ with respect to $(M,\eta)$ is 
$$ent_{i}(\overline{\sigma \vee \tau})_{(M,\eta)}=sup\Big \{H_{i}(M',\eta) \mid M'\in \mathcal{F}_{i}\big( (\sigma \vee \tau)(M)\big) \Big \},$$
that is, where $M'\leq \sigma(M) + \tau(M)$ and $i(M')<\infty$.

\item The $i$-entropy for $(\overline{\sigma \cdot \tau})$ with respect to $(M,\eta)$ is 
$$ent_{i}(\overline{\sigma \cdot \tau})_{(M,\eta)}=sup \Big \{H_{i}(M',\eta) \mid M'\in \mathcal{F}_{i}\big( (\sigma \cdot \tau)(M)\big) \Big \},$$
that is, $M'\leq \sigma(\tau(M))$ and where $i(M')<\infty$.
    
\item The $i$-entropy for $(\overline{\sigma : \tau})$ with respect to $(M,\eta)$ is  
$$ent_{i}(\overline{\sigma : \tau})_{(M,\eta)}=sup\Big \{H_{i}(M',\eta) \mid M'\in \mathcal{F}_{i}\big( (\sigma : \tau)(M)\big) \Big \},$$
that is, where $M'\leq \tau(M/\sigma(M))$ and $i(M')<\infty$.
\end{itemize}
\end{prop}

\begin{prop}\label{Prop entropy} 
Let $\sigma, \tau \in R\textit{-}pr$ and let $(M,\eta) \in \textbf{Flow}_{R\textit{-}Mod}$. Then, we have that  
\begin{center}
    $ent_{i}(\overline{\sigma \cdot \tau})_{(M,\eta)}\leq_{\mathbb{R}} ent_{i}(\overline{\sigma \wedge \tau})_{(M,\eta)} \leq_{\mathbb{R}} ent_{i}(\overline{\sigma \vee \tau})_{(M,\eta)} \leq_{\mathbb{R}} ent_{i}(\overline{\sigma: \tau})_{(M,\eta)}$.
\end{center}
where $\overline{\sigma}$ represents the induced preradical by $\sigma$, via the assignment $\phi$. 
\end{prop}

\begin{dem}
Let $\sigma, \tau\in R\textit{-}pr$ and let $(M,\eta) \in \textbf{Flow}_{R\textit{-}Mod}$.  
\begin{itemize}
\item[(i)] 
Since $(\sigma\cdot \tau)(M)\leq (\sigma \wedge \tau)(M)$ for all $M\in R\textit{-}Mod$, it follows that 
\begin{center}
    $\mathcal{F}_{i}\big( (\sigma \cdot \tau)(M)\big)\subseteq \mathcal{F}_{i}\big( (\sigma \wedge \tau)(M)\big)$.
\end{center}        
    The latter in turn implies that 
\begin{center}
    $\Big \{H_{i}(M',\eta) \mid M'\in \mathcal{F}_{i}\big( (\sigma \cdot \tau)(M)\big) \Big \}\leq_{\mathbb{R}} \Big\{H_{i}(M',\eta) \mid M'\in \mathcal{F}_{i}\big( (\sigma \wedge \tau)(M)\big) \Big \}$.
\end{center}
Therefore, 
\begin{center}
     $sup \Big \{H_{i}(M',\eta) \mid M'\in \mathcal{F}_{i}\big( (\sigma \cdot \tau)(M)\big) \Big \}\leq_{\mathbb{R}} sup \Big\{H_{i}(M',\eta) \mid M'\in \mathcal{F}_{i}\big( (\sigma \wedge \tau)(M)\big) \Big \}$.
\end{center}
Hence, $ent_{i}(\overline{\sigma \cdot \tau})_{(M,\eta)}\leq_{\mathbb{R}} ent_{i}(\overline{\sigma \wedge \tau})_{(M,\eta)}$.     
    
\item[(ii)]  
Since $(\sigma \wedge \tau)(M)\leq (\sigma\vee \tau)(M)$ for all $M\in R\textit{-}Mod$, then we have that  
\begin{center}
    $\mathcal{F}_{i}\big( (\sigma \wedge \tau)(M)\big)\subseteq \mathcal{F}_{i}\big( (\sigma \vee \tau)(M)\big)$.    
\end{center}    
Thus,   
\begin{center}
    $\Big \{H_{i}(M',\eta) \mid M'\in \mathcal{F}_{i}\big( (\sigma \wedge \tau)(M)\big) \Big \} \leq_{\mathbb{R}} \Big\{H_{i}(M',\eta) \mid M'\in \mathcal{F}_{i}\big( (\sigma \vee \tau)(M)\big) \Big \}$, 
\end{center}   
which implies that 
\begin{center}
     $sup \Big \{H_{i}(M',\eta) \mid M'\in \mathcal{F}_{i}\big( (\sigma \wedge \tau)(M)\big) \Big \}\leq_{\mathbb{R}} sup \Big\{H_{i}(M',\eta) \mid M'\in \mathcal{F}_{i}\big( (\sigma \vee \tau)(M)\big) \Big \}$.
\end{center}
Therefore, $ent_{i}(\overline{\sigma \wedge \tau})_{(M,\eta)} \leq_{\mathbb{R}} ent_{i}(\overline{\sigma \vee \tau})_{(M,\eta)}$
    
\item[(iii)] Since $(\sigma\vee \tau)(M)\leq (\sigma:\tau)(M)$, for all $M\in R\textit{-}Mod$, then 
\begin{center}
    $\mathcal{F}_{i}\big( (\sigma \vee  \tau)(M)\big)\subseteq \mathcal{F}_{i}\big( (\sigma : \tau)(M)\big)$.    
\end{center} 
Thus, 
\begin{center}
    $\Big \{H_{i}(M',\eta) \mid M'\in \mathcal{F}_{i}\big( (\sigma \vee \tau)(M)\big) \Big \}\leq_{\mathbb{R}} \Big\{H_{i}(M',\eta) \mid M'\in \mathcal{F}_{i}\big( (\sigma : \tau)(M)\big) \Big \}$,
\end{center}   
which in turn implies 
\begin{center}
     $sup \Big \{H_{i}(M',\eta) \mid M'\in \mathcal{F}_{i}\big( (\sigma \vee \tau)(M)\big) \Big \}\leq_{\mathbb{R}} sup \Big\{H_{i}(M',\eta) \mid M'\in \mathcal{F}_{i}\big( (\sigma : \tau)(M)\big) \Big \}$.
\end{center}
Hence, $ent_{i}(\overline{\sigma \vee \tau})_{(M,\eta)} \leq_{\mathbb{R}} ent_{i}(\overline{\sigma: \tau})_{(M,\eta)}$.
\end{itemize}
\qed
\end{dem}

We close this subsection with a result that relates the algebraic entropy $log$ defined on categories of the form $S$-Mod and $R$-Mod, whenever there is a ring homomorphism between the base rings. With this in mind, we will show how the entropy values for preradicals on $S$-Mod relate to the entropy values of preradicals on $R$-Mod, whenever the ring homomorphism $t:R\longrightarrow S$ is surjective. 
\\ Let us then suppose that $t\colon R\longrightarrow S$ denotes a ring homomorphism. Then each $S$-module $M$ admits an $R$-module structure via the homomorphism $t$. Indeed, given $M\in S$-Mod we have that $t$ induces an action on $M$ - as an $R$-module-, given by the correspondence rule $r\cdot m:=t(r)\cdot m$, for any $r\in R$ and $m\in M$. Moreover, $t$ induces a functor $F_{t}:S\textit{-}Mod \longrightarrow R\textit{-}Mod$ that assigns to each $M\in S\textit{-}Mod$, the $R$-module $F_{t}(M)=M$, and assigns to each morphism $f:M\longrightarrow M'$ in $S\textit{-}Mod$, the $R$-morphism $F_{t}(f)=f$. 

\begin{prop}
Let $log_{s}$ and $log_{r}$ be the invariant $log()$ defined on $S\textit{-}Mod$ and $R\textit{-}Mod$ respectively. If $t\colon R\longrightarrow S$ is a ring homomorphism, then for any $M\in S\textit{-}Mod$ and its respective $R$-module $F_{t}(M)$, one has that  
\begin{center}
$ent_{log_{s}}(M)\leq ent_{log_{r}}(F_{t}(M))$. 
\end{center}
\end{prop}

\begin{dem}
Let $t\colon R\longrightarrow S$ be a ring homomorphism. For the sake of clarity, we will denote an $S$-module by $M_{S}$ and to its image under functor $F_{t}$ by $M_{R}$. We first note that, if $L_{S}\in \mathcal{F}_{log_{s}}(M_{S})$ then $L_{R}\in \mathcal{F}_{log_{r}}(M_{R})$. Indeed, the functor $F_{t}$ send the inclusion map $\iota_{S}:L_{S}\longrightarrow M_{S}$ into the inclusion map $i_{R}:L_{R}\longrightarrow M_{R}$ in $R$-Mod. Hence, as $L_{S}\in \mathcal{F}_{log_{s}}(M_{S})$ then $log(L_{R})=Log|L|=log(L_{S})<  \infty$, and thus $L_{R}\in \mathcal{F}_{log_{r}}(M_{R})$ for each $L_{S}\in \mathcal{F}_{log_{s}}(M_{S})$. Consequently, we have that  for each endomorphism $\eta:M_{S}\longrightarrow M_{S}$ and each $L_{S}\in \mathcal{F}_{log_{s}}(M_{S})$ one has that 
\begin{align*}
T_{n}(L_{S},\eta)=L_{S}+\eta(L)_{S}+\cdots + \eta^{n-1}(L)_{S} \\ 
= L_{R}+ \eta(L)_{R} + \cdots +  \eta^{n-1}(L)_{R}= T_{n}(L_{R},\eta_{R}),
\end{align*}
where $\eta_{R}$ denotes the endomorphism $\eta_{R}:M_{R}\longrightarrow M_{R}$ resulting from evaluating functor $\mathcal{F}_{t}$ on $\eta$. Considering the above we have that 
\begin{center}
$H_{log_{s}}(L_{S},\eta)= \lim_{n\to \infty} \frac{log|(T_{n}(L_{S},\eta))|}{n}$ \vspace{2mm}
$=\lim_{n\to \infty} \frac{log|(T_{n}(L_{R},\eta_{R}))|}{n}=H_{log_{r}}(L_{R},\eta_{R})$,
\end{center}
and thus,  
\begin{center}
$ent_{log_{s}}(\eta) = sup \{ H_{log_{s}}(L_{S},\eta) \mid L_{S}\in \mathcal{F}_{log_{s}}(M_{S}) \}$ \\ 
\vspace{3mm}
$=sup \{ H_{log_{r}}(L_{R},\eta_{R}) \mid  L_{R}\in \mathcal{F}_{log_{r}}(M_{R}) \}$ \\ 
\vspace{3mm}
$\leq_{\mathbb{R}} sup \{ H_{log_{r}}(N,\eta_{R}) \mid N\in \mathcal{F}_{log_{r}}(M_{R}) \}= ent_{log_{r}}(\eta_{R}).$
\end{center}

Therefore,
\begin{center}
$ent_{log_{s}}(M_{S}) = sup \{ent_{log_{s}}(\eta) \mid \eta \in End_{S} (M_{S}) \}$ \\
\vspace{3mm}
$\leq_{\mathbb{R}} sup \{ent_{log_{r}}(\eta_{R}) \mid \eta_{R} \in End_{R} (M_{R}) \}$ \\
\vspace{3mm}    
$\leq_{\mathbb{R}} sup \{ ent_{log_{r}}(\mu) \mid \mu \in End_{R} (M_{R}) \}= ent_{log_{r}}(M_{R}).$
\end{center}
% Observe that the equality holds when $t$ is an epimorphism, since the functor $\mathcal{F}_{t}$ induces a surjective function from $Hom(M_{S},M_{S})=End_{S}(M_{S})$ to $Hom(M_{R},M_{R})=End_{R}(M_{R})$.
\qed 
\end{dem}

In case that $t\colon R\longrightarrow S$ is a surjective ring homomorphism, one always has that the induced functor $F_{t}\colon S\textit{-}Mod \longrightarrow R\textit{-}Mod$ is full (\cite[Theorem 2.5]{Pre 2}). Furthermore, as $F_{t}$ respects inclusion maps and is injective on objects, by \cite[Theorem 2.4]{Pre 2}  we have that $F_{t}$ induces an injective assignment $\phi_{t}:S\textit{-}pr\longrightarrow R\textit{-}pr$ that is order preserving. This assignment is given by
\[ 
\phi_{t}(\sigma) = \underset{M\in S\textit{-}Mod}{\bigvee} \  \alpha_{\sigma(M)}^{M}  
\]
where $\alpha_{\sigma(M)}^{M}$ denotes the alpha preradical on $R$-pr defined by the objects $\sigma(M)=F_{t}(\sigma(M))$ and $M=F_{t}(M)$. With these facts in mind, we have the next

\begin{prop}
Let $log_{s}$ and $log_{r}$ be the invariant $log()$ defined on $S\textit{-}Mod$ and $R\textit{-}Mod$ respectively. If $t\colon R\longrightarrow S$ is a ring epimorphism, then 
\begin{center}
$ent_{log_{s}}(\sigma)|_{M_{S}} \leq_{\mathbb{R}} ent_{log_{r}}(\phi_{t}(\sigma))|_{M_{R}}$    
\end{center}
for each preradical $\sigma \in S\textit{-}pr$. 
\\ Here, $M_{S}\in S\textit{-}Mod$ and $M_{R}$ is the image of $M_{S}$ under functor $F_{t}$.
\end{prop}

\begin{dem}
Let $M_{S}\in S\textit{-}Mod$ and let  $M_{R}=F_{t}(M_{S})$. 
In order to avoid confusion, we will denote $\sigma(M_{S})=\sigma(M)_{S}$ and thus, we denote  $\sigma(M)_{R}$ to its corresponding $R$-module under functor $F_{t}$. Now, for $\sigma(M)_{R}\leq M_{R}$ one has that 
\begin{center}
$\alpha_{\sigma(M_{R})}^{M_{R}} \leq \underset{M\in S\textit{-}Mod}{\bigvee} \  \alpha_{\sigma(M)}^{M}$. \end{center}
Further, we have that $\alpha_{\sigma(M)_{R}}^{M_{R}}(M_{R})=\sigma(M)_{R}$. Therefore, since $M_{S}'\in \mathcal{F}_{log_{s}}(\sigma(M_{S}))$ implies that $M_{R}'\in \mathcal{F}_{log_{r}}(\alpha_{\sigma(M_{R})}^{M_{R}}(M_{R}))\subseteq \mathcal{F}_{log_{r}}\Big( \big( \underset{M\in S\textit{-}Mod}{\bigvee} \  \alpha_{\sigma(M)}^{M}\big)(M_{R}) \Big)$, we then have that  
\[
ent_{log_{s}}(\sigma)|_{M_{S}} \leq_{\mathbb{R}} ent_{log_{r}}(\alpha_{\sigma(M_{R})}^{M_{R}})|_{M_{R}} \leq_{\mathbb{R}} ent_{log_{r}}(\phi_{t}(\sigma))|_{M_{R}}
\]

\qed
\end{dem}

\subsection{Functorial entropy between $\mathcal{L}\textit{-}pr$ and $\mathcal{L_{M}}\textit{-}pr$}

We will now use the notion of functorial entropy described in \cite{Norm} to define entropy for \textit{lattice preradicals} (See \cite{Albu 2}). In doing this, we will show some results that relates entropy for lattice preradicals with entropy for module preradicals. 

It is well known that each $R$-module $M$ induces a complete modular lattice $L(M)$, whose elements correspond to all submodules $N$ of $M$. The meet and the join operations in $L(M)$ ( denoted by $\wedge$ and $\vee$ respectively) are given by the intersection and the sum of submodules. Further, each $L(M)$ is bounded above by $M$ and below by $\{0\}$. As the authors show in \cite{Albu 1}, the collection of all bounded modular lattices are the objects of a category, denoted by $\mathcal{L_{M}}$, whose morphisms are called \textit{linear morphisms}: these summon the property of having a kernel $ker(\eta)$ for every morphisms $\eta:M\longrightarrow M'$ in $R$-mod, so that the First Isomorphisms Theorem holds; this is:
\begin{center}
$M/ker(\eta) \cong Im(\eta)$.      
\end{center}
Moreover, for each $R$-morphisms $M \overset{\eta}{\longrightarrow} M'$ one can define a \textit{linear morphism} $L(M)\overset{\mathcal{X}_{\eta}}{\longrightarrow}L(M')$ as follows: if $N\in L(M)$ -that is, $N$ is a submodule of $M$- then $\mathcal{X}_{\eta}(N)=\eta(N)$, where $\eta(N)$ is the submodule of $M'$ generated by the image of $N$ under $\eta$. The above can be summarized by saying that, for each unital ring $R$, one has a functor $\mathcal{X}:R\textit{-}Mod\longrightarrow \mathcal{L_{M}}$ with correspondence $M\longmapsto L(M)$ on objects and $\eta \longmapsto \mathcal{X}_{\eta}$ on morphisms. We note here that the image of the functor $\mathcal{X}$ defines a non-full subcategory of $\mathcal{L_{M}}$, which we will denote by $\mathcal{SL_{M}}$. In this case, the set of morphisms from $L(M)$ to $L(M')$ within $\mathcal{SL_{M}}$ is given by the set 
\begin{center}
$Hom_{\mathcal{SL_{M}}}\big (L(M),L(M') \big )= \Big \{ \mathcal{X}_{\eta} \mid M\overset{\eta}{\longrightarrow} M' \mbox{ is a morphism in }R\textit{-}Mod \Big \}$.    
\end{center}  

When considering preradicals, in \cite{Albu 2} the authors show the existence of an assignment, denoted here by $\phi:\mathcal{L_{M}}\textit{-}pr\longrightarrow R\textit{-}pr$, between the big lattice of lattice preradicals and the big lattice of $R$-module preradicals. This assignment is given by the correspondence $\Tilde{\sigma}\longmapsto \sigma$, where $\sigma$ evaluated in an $R$-module $M$ coincides with the respective submodule obtained by evaluating $\Tilde{\sigma}$ in $L(M)$. In other words, for any $M\in R$-Mod, we have $\sigma(M)=X_{L(M)}^{\Tilde{\sigma}}$ where $X_{L(M)}^{\Tilde{\sigma}}$ is such that  $\Tilde{\sigma}(L(M))=X_{L(M)}^{\Tilde{\sigma}}/0$\footnote{A lattice preradical on $\mathcal{L_{M}}$ assigns to each lattice $L$ a sublattice of the form $X/0$.  Here, we have written $X_{L(M)}^{\Tilde{\sigma}}$ with a capital letter to emphasize that $X_{L(M)}^{\Tilde{\sigma}}$ defines a submodule of $M$. }. Lastly, as shown in \cite{Sebas}, $\phi$ is an order preserving assignment which also respects the join, the product and the coproduct operations. 

Let us now suppose that $i$ is a subadditive invariant on $R$-Mod, and let $M\in R$-Mod. As noted before, $M$ induces a complete modular lattice $\big( L(M), \wedge, \vee, \subseteq \big)$. This complete modular lattice contains the set 
\begin{center}
$\mathcal{F}_{i}(L(M))=\{N\in L(M) \mid i(N)<\infty \}$.    
\end{center}
Further, as $i$ is an subadditive invariant, one has that 
\begin{center}
$i(N_{1}+N_{2})\leq i(N_{1}) + i(N_{2})< \infty$.    
\end{center} 
Hence, $N_{1}\vee N_{2}=N_{1}+N_{2} \in \mathcal{F}_{i}(L(M))$, and thus the tuple $(\mathcal{F}_{i}(L(M)), \vee, \subseteq)$ defines a \textit{join semi-lattice} \footnote{A semi-lattice is a set $S$ with a binary operation $*$ which is commutative, associative, and satisfies $s*s=s$ for all $s\in S$. In particular, any lattice $L$ is a semi-lattice under $\vee$ and also under $\wedge$.}. Moreover, according to \cite[Section 5]{Norm}, we have that the tuple $(\mathcal{F}_{i}(L(M)), \vee, \subseteq, i)$ defines a subadditive \textit{normed semi-lattice}\footnote{ A normed semi-lattice is a semi-lattice $S$ equipped with a norm $i:S\longrightarrow \mathbb{R}_{\geq 0}$}. 
\begin{prop}
Let $i$ be an invariant on $R$-Mod, and let $\eta\colon M\longrightarrow M'$ be an $R$-morphism. Then $\eta$ induces a morphism $\mathcal{F}_{i}(\eta)$ between the join semi-lattices $\mathcal{F}_{i}(L(M))$ and $\mathcal{F}_{i}(L(M'))$. This morphism is given by $\big(
 \mathcal{F}_{i}(\eta)\big)(N)=\eta(N)$ for each $N\in \mathcal{F}_{i}(L(M))$.  
\end{prop}
\begin{dem}
Let $i$ be an invariant on $R$-Mod, and let $\eta\colon M\longrightarrow M'$ be an $R$-morphism. We will first see that $\mathcal{F}_{i}(\eta)\colon \mathcal{F}_{i}(L(M))\longrightarrow \mathcal{F}_{i}(L(M'))$, with correspondence rule $\big(
 \mathcal{F}_{i}(\eta)\big)(N)=\eta(N)$ for each $N\in \mathcal{F}_{i}(L(M))$, is well defined. Indeed, as $i$ is an invariant, for each $N\in \mathcal{F}_{i}(L(M))$ we have that 
\[
i(\eta(N))=i(N/ker(\eta))\leq i(N)<\infty. 
\]
Thus, $i(\eta(N))<\infty$, which implies that  $\big(
 \mathcal{F}_{i}(\eta)\big)(N)=\eta(N)\in \mathcal{F}_{i}(L(M'))$. 
\\ Now, if $N_{1},N_{2}\in \mathcal{F}_{i}(L(M))$, then
\begin{align*}
\big( \mathcal{F}_{i}(\eta)\big)(N_{1}\vee N_{2}) = \eta(N_{1}\vee N_{2})=\eta(N_{1} + N_{2})
=\eta(N_{1})+\eta(N_{2})
\\ =\big( \mathcal{F}_{i}(\eta)\big)(N_{1}) +\big( \mathcal{F}_{i}(\eta)\big)(N_{2})=\big( \mathcal{F}_{i}(\eta)\big)(N_{1}) \vee \big( \mathcal{F}_{i}(\eta)\big)(N_{2})\in \mathcal{F}_{i}(L(M')).
\end{align*}
Therefore, $\mathcal{F}_{i}(\eta):\mathcal{F}_{i}(L(M))\longrightarrow \mathcal{F}_{i}(L(M'))$ is a morphism.
\qed 
\end{dem}

In this way, we have a functor $\mathcal{S}\colon R\textit{-Mod}\longrightarrow \mathcal{L}^{\dag}$, where $\mathcal{L}^{\dag}$ is the category whose objects are normed join semi-lattices and whose morphisms are the contractive\footnote{These morphisms are named contractive since  $i(\eta(N))<i(N)$, for all $N$.} join semi-lattice morphisms. Here $\mathcal{S}$ assigns to each $R$-module $M$ the join semi-lattice $(\mathcal{F}_{i}(L(M)), \vee, \subseteq, i)$, and assigns to each $R$-morphism $\eta$ the join semi-lattice morphism $\mathcal{F}(\eta)$. Similarly, we can define a functor $\mathcal{F}\colon \mathcal{SL_{M}}\longrightarrow \mathcal{L}^{\dag}$ where $\mathcal{F}$ assigns to each object $\big( L(M), \wedge, \vee, \subseteq \big)$ the object $(\mathcal{F}_{i}(L(M)), \vee, \subseteq, i)$ and assigns to each morphism $\mathcal{X}_{\eta}$ the join semi-lattice morphism $\mathcal{F}(\eta)$. As we now see, the entropy for objects and endomorphisms within the category $\mathcal{SL_{M}}$ is equivalent to that defined on objects and endomorphisms within $R$-Mod:

\begin{defi} \label{entropy chi}
Let $i$ be a subadditive invariant on $R$-Mod, $M$ an $R$-module and $L(M)$ the complete lattice of all submodules of $M$. If $\eta:M\longrightarrow M$ is an endomorphism and $\mathcal{X}_{\eta}$ is the induced endomorphism of $L(M)\in \mathcal{SL_{M}}$, then the algebraic $i$-entropy of $\mathcal{X}_{\eta}$ with respect to $N\in \mathcal{F}_{i}(L(M))$ is 
\begin{center}
$H_{i}(N,\mathcal{X}_{\eta})= \lim_{n\to \infty} \frac{i(T_{n}(N,\mathcal{X}_{\eta}))}{n}$.
\end{center} 
Hence, the algebraic $i$-entropy of $\mathcal{X}_{\eta}$ is 
\begin{center}
    $ent_{i}(\mathcal{X}_{\eta})=sup\Big \{H_{i}(N,\mathcal{X}_{\eta}) \mid N\in \mathcal{F}_{i}(L(M))\Big \}$.
\end{center}
\end{defi}

A first consequence of Definition \ref{entropy chi} is that the entropy for a complete modular lattice of the form $L(M)\in \mathcal{SL_{M}}$ (with $M\in R$-Mod) is:  
\begin{center}
$ent_{i}(L(M))=sup \big \{ ent_{i}(\mathcal{X}_{\eta}) \mid \mathcal{X}_{\eta} \in End_{\mathcal{SL_{M}}}(L(M)) \big \}$ \vspace{+3mm}
\\ $= sup \big \{ ent_{i}(\mathcal{X}_{\eta}) \mid \eta \in End(M) \big \} $.
\end{center}
Therefore, as $N\in \mathcal{F}_{i}(L(M))$ if, and only if, $N\in \mathcal{F}_{i}(M)$ and the fact that $\mathcal{X}_{\eta}(N)=\eta(N)$, it follows the next 
\begin{prop}
Let $M\in R\textit{-}Mod$ and $L(M)\in \mathcal{L_{M}}$. Then 
\begin{center}
$ent_{i}(M)= ent_{i}(L(M))$.    
\end{center}
\end{prop}

Furthermore, according to \cite[Theorem 5.26]{Norm}, we have that for each endomorphism $M\overset{\eta}{\longrightarrow} M$ in $R$-Mod,
\begin{center}
$\textbf{h}_{\mathcal{S}}(\eta)=ent_{i}(\eta)=ent_{i}(\mathcal{X}_{\eta})=\textbf{h}_{\mathcal{F}}$,
\end{center}
where $\textbf{h}_{\mathcal{S}}$ is the functorial entropy defined by the functor $\mathcal{S}:R\textit{-}Mod \longrightarrow \mathcal{L}^{\dag}$ and $\textbf{h}_{\mathcal{F}}$ is the functorial entropy defined by the functor $\mathcal{F}:\mathcal{SL_{M}}\longrightarrow \mathcal{L}^{\dag}$. 

Considering the above, we will define entropy for lattice preradicals on objects of the category $\mathcal{SL_{M}}$, and we will show that this definition is compatible with the definition of entropy for preradicals on $R$-Mod.  

\begin{defi}
Let $i$ be a subadditive invariant on $R$-Mod, $M$ an $R$-module, $L(M)$ the complete lattice of all submodules of $M$, and $\Tilde{\sigma}$ a lattice preradical on $\mathcal{L_{M}}$. Then the algebraic $i$-entropy of $\Tilde{\sigma}$ with respect to $L(M)\in \mathcal{SL_{M}}$ is
\[
ent_{i}(\Tilde{\sigma})|_{L(M)}\colon =\underset{\mathcal{X}_{\eta}\in End_{\mathcal{SL_{M}}}(L(M))}{sup} \Big \{ H_{i}(N,\mathcal{X}_{\eta})  \mid N \in \mathcal{F}_{i}(\Tilde{\sigma}(L(M)))  \Big \}
\]
where $H_{i}(N,\mathcal{X}_{\eta})=\lim_{n\to \infty} \frac{i(T_{n}(N,\mathcal{X}_{\eta}))}{n}$.
\end{defi}

\begin{prop} \label{ent equality}
Let $i$ be an invariant on $R$-Mod, and let $\phi\colon \mathcal{L_{M}}\textit{-}pr\longrightarrow R\textit{-}pr$. If $\Tilde{\sigma}\in \mathcal{L_{M}}\textit{-}pr$ and $\phi(\Tilde{\sigma})=\sigma$ denotes its induced preradical on $R$-pr, then given $M\in R\textit{-}Mod$ and $L(M)\in \mathcal{SL_{M}}$ one has that
\begin{center}
$ent_{i}(\sigma)|_{M}= ent_{i}(\Tilde{\sigma})|_{L(M)}$.
\end{center}
\end{prop}

\begin{dem}
First observe that $N\in \mathcal{F}_{i}(\sigma(M))$ if, and only if, $N\in \mathcal{F}_{i}(\Tilde{\sigma}(L(M)))$. Indeed, by definition of $N\in \mathcal{F}_{i}(\sigma(M))$ one has that $N$ is a submodule of $\sigma(M)=X_{L(M)}^{\Tilde{\sigma}}$ with $i(N)<\infty$. Recall that $X_{L(M)}^{\Tilde{\sigma}}$ is such that 
\begin{center}
$\Tilde{\sigma}(L(M))=\Big( X_{L(M)}^{\Tilde{\sigma}}/0 \Big)$.    
\end{center} 
This in turn implies that $N$ is an element of the lattice $\Big( X_{L(M)}^{\Tilde{\sigma}}/0 \Big)$, and hence $N\in \mathcal{F}_{i}(\Tilde{\sigma}(L(M)))$. 
\\ Conversely, if $N\in \mathcal{F}_{i}(\Tilde{\sigma}(L(M)))$ then $i(N)<\infty$ and $N$ is an element of the lattice $\Big( X_{L(M)}^{\Tilde{\sigma}}/0 \Big)$. Since by definition $X_{L(M)}^{\Tilde{\sigma}} = \sigma(M)$, then $N\leq \sigma(M)$ and thus, $N\in \mathcal{F}_{i}(\sigma(M))$.

Lastly, since  
\[
H_{i}(N,\eta)=\lim_{n\to \infty} \frac{i(T_{n}(N,\eta))}{n} =\lim_{n\to \infty} \frac{i(T_{n}(N,\mathcal{X}_{\eta}))}{n}= H_{i}(N,\mathcal{X}_{\eta}),
\]
for each endomorphism $M\overset{\eta}{\longrightarrow} M$ in $R$-Mod and its induced endomorphism $L(M)\overset{\mathcal{X}_{\eta}}{\longrightarrow} L(M)$ in $\mathcal{SL_{M}}$, we have that  
\begin{center}
$ent_{i}(\sigma)|_{M} =\underset{\eta\in End(M)}{sup} \Big \{ H_{i}((N,\eta)  \mid N\in \mathcal{F}_{i}(\sigma(M)) \Big \}$
\\ \vspace{2mm}
$= \underset{\mathcal{X}_{\eta}\in End_{\mathcal{SL_{M}}}(L(M))}{sup} \Big \{ H_{i}(N,\mathcal{X}_{\eta})  \mid N \in \mathcal{F}_{i}(\Tilde{\sigma}(L(M)))  \Big \}$
\\ \vspace{2mm}
$= ent_{i}(\Tilde{\sigma})|_{L(M)}$.
\end{center}
\qed 
\end{dem}

Finally, since the assignment $\phi\colon \mathcal{L_{M}}\textit{-}pr\longrightarrow R\textit{-}pr$ preserves the four operations between preradicals (see \cite[Section 8]{Sebas}), it is straightforward the next

\begin{coro} \label{Equalities}
Let $i$ be an invariant on $R$-Mod, and let $\Tilde{\sigma}\in \mathcal{L_{M}}\textit{-}pr$ with $\phi(\Tilde{\sigma})=\sigma\in R\textit{-}pr$ its induced preradical under the assignment $\phi$. If $M\in R\textit{-}Mod$ and $L(M)\in \mathcal{SL_{M}}$ then  
\begin{itemize}
\item[(i)]  $ent_{i} (\Tilde{\sigma} \wedge \Tilde{\tau})|_{L(M)} = ent_{i} (\sigma \wedge \tau)|_{M}$ 
\item[(ii)]  $ent_{i} (\Tilde{\sigma} \vee \Tilde{\tau})|_{L(M)}= ent_{i} (\sigma \vee \tau)|_{M}$
\item[(iii)] $ent_{i} (\Tilde{\sigma} \cdot \Tilde{\tau})|_{L(M)}= ent_{i} (\sigma \cdot \tau)|_{M}$
\item[(iv)] $ent_{i} (\Tilde{\sigma} : \Tilde{\tau})|_{L(M)}= ent_{i} (\sigma : \tau)|_{M}$
\end{itemize}
\end{coro}

\begin{coro}
Let $i$ be an invariant on $R$-Mod, and let $\Tilde{\sigma}\in \mathcal{L_{M}}\textit{-}pr$ with $\phi(\Tilde{\sigma})=\sigma\in R\textit{-}pr$ its induced preradical under the assignment $\phi$. If $M\in R\textit{-}Mod$ and $L(M)\in \mathcal{SL_{M}}$ then  
\begin{center}
    $ent_{i}(\Tilde{\sigma}\cdot \Tilde{\tau})|_{L(M)}\leq_{\mathbb{R}} 
    ent_{i}(\Tilde{\sigma}\wedge \Tilde{\tau})|_{L(M)} \leq_{\mathbb{R}}
    ent_{i}(\Tilde{\sigma}\vee \Tilde{\tau})|_{L(M)} \leq_{\mathbb{R}} 
    ent_{i}(\Tilde{\sigma} : \Tilde{\tau})|_{L(M)}$.
\end{center}
\end{coro}

\begin{dem}
It follows from Corollary \ref{Inequality} and Corollary \ref{Equalities}
\qed
\end{dem}

\section*{Acknowledgments}
This work was supported by unrestricted funds to the Center for Engineered Natural Intelligence at the University of California San Diego.

\end{document}